\documentclass{gtart_h}

\def\ifplaintex{\expandafter\ifx\csname documentclass\endcsname\relax}

\def\gtp{{\mathsurround=0pt\it $\cal G\mskip-2mu$eometry \&\ 
$\cal T\!\!$opology $\cal P\!$ublications}}  

\def\Addressesr{\bigskip
{\small \parskip 0pt \leftskip 0pt \rightskip 0pt plus 1fil \def\\{\par}
\sl\theaddress\par
\medskip
\rm Email:\stdspace\tt\theemail\hfill\rm Received:\qua\receiveddate \par}}

\def\recd{{\small Received:\qua\receiveddate\ifx\reviseddate\relax
\else\qquad Revised:\qua\reviseddate\fi\par}} 

\def\lognumber#1{\def\thelognumber{#1}}
\def\volumenumber#1{\def\thevolumenumber{#1}}
\def\volumeyear#1{\def\thevolumeyear{#1}}
\def\papernumber#1{\def\thepapernumber{#1}}
\def\pagenumbers#1#2{\def\startpage{#1}\def\finishpage{#2}}
\def\published#1{\def\publishdate{#1}}

\def\received#1{\def\receiveddate{#1}}

\def\accepted#1{\def\accepteddate{#1}}
\def\asciititle#1{\def\theasciititle{#1}}

\long\def\asciiabstract#1{\long\def\theasciiabstract{#1}}

\let\\\par\let\thelognumber\relax\let\thevolumenumber\relax
\let\thepapernumber\relax\let\thevolumeyear\relax\let\startpage\relax
\let\finishpage\relax\let\publishdate\relax\let\receiveddate\relax
\let\reviseddate\relax\let\accepteddate\relax\let\theasciititle\relax
\let\theasciiauthors\relax
\let\theasciiabstract\relax

\let\theasciiemail\relax

\ifplaintex
\font\logobig=cmssbx10 scaled 3836
\font\logomed=cmssbx10 scaled 2557
\else
\font\logobig=cmssbx10 scaled 4200
\font\logomed=cmssbx10 scaled 2800
\fi

\long\def\makeagttitle{   
\count0=\startpage
\agt\hfill      
\hbox to 45truept{\vbox to 0pt{\vglue -13truept{\logomed A\kern -.37em{\logobig 
T}\kern -.38em G}\vss}\hss}
\break
{\small Volume \thevolumenumber\ (\thevolumeyear)
\startpage--\finishpage\nl
Published: \publishdate}

\vglue .25truein

{\parskip=0pt\leftskip 0pt plus
1fil\def\\{\par\smallskip}{\Large\bf\thetitle}\par\medskip} \vglue
0.05truein

{\parskip=0pt\leftskip 0pt plus 1fil\def\\{\par}{\sc\theauthors}
\par\medskip}%
 
\vglue 0.03truein 

{\small\leftskip 25truept\rightskip 25truept{\bf Abstract}\stdspace\theabstract

{\bf AMS Classification}\stdspace\theprimaryclass
\ifx\thesecondaryclass\relax\else; \thesecondaryclass\fi\par
{\bf Keywords}\stdspace \thekeywords\par}\vglue 7truept

}   

\ifplaintex
\font\phead=cmsl9 scaled 950
\font\pnum=cmbx10 scaled 913
\font\pfoot=cmsl9 scaled 950
\headline{\vbox to 0pt{\vskip -4.5mm\line{\small\phead\ifnum
\count0=\startpage ISSN 1472-2739 (on-line) 1472-2747 (printed)
\hfill {\pnum\folio}\else\ifodd\count0\def\\{ }%
\ifx\theshorttitle\relax\thetitle\else\theshorttitle\fi\hfill{\pnum\folio}
\else\def\\{ and }{\pnum\folio}\hfill\ifx\theshortauthors\relax\theauthors
\else\theshortauthors\fi\fi\fi}\vss}}
\footline{\vbox to 0pt{\vglue 0mm\line{\small\pfoot\ifnum\count0=\startpage
\copyright\ \gtp\hfill\else
\agt, Volume \thevolumenumber\ (\thevolumeyear)\hfill\fi}\vss}}
\else
\font=cmsl9 scaled 1050
\font\lnum=cmbx10 
\font=cmsl9 scaled 1050
\makeatletter
\def\@oddhead{{\small\lhead\ifnum\count0=\startpage ISSN 1472-2739 
(on-line) 1472-2747 (printed)\hfill {\lnum\number\count0}\else\ifodd\count0
\def\\{ }\ifx\theshorttitle\relax \thetitle \else\theshorttitle\fi\hfill
{\lnum\number\count0}\else\def\\{ and }{\lnum\number\count0}
\hfill\ifx\theshortauthors\relax 
\theauthors\else\theshortauthors\fi\fi\fi}}\def\@evenhead{\@oddhead}
\def\@oddfoot{\small\lfoot\ifnum\count0=\startpage\copyright\ \gtp\hfill\else
\agt, Volume \thevolumenumber\ (\thevolumeyear)\hfill\fi}
\def\@evenfoot{\@oddfoot}
\makeatother
\fi
\let\maketitlepage\makeagttitle
\let\makeshorttitle\maketitlepage
\let\maketitle\maketitlepage

\newwrite\gtoutfile
\long\gdef\makeheadfile{  
{\def\\{, }\def\s{ }
\immediate\openout\gtoutfile head.xxx
\immediate\write\gtoutfile{Proxy-for: \ifx\theasciiauthors\relax
\theauthors\else\theasciiauthors\fi\s<\ifx\theasciiemail\relax\theemail\else\theasciiemail\fi>}
\immediate\write\gtoutfile{\noexpand\\}
\immediate\write\gtoutfile{Authors: \ifx\theasciiauthors\relax
\theauthors\else\theasciiauthors\fi}
{\def\\{ }\immediate\write\gtoutfile{Title: \ifx\theasciititle\relax
\thetitle\else\theasciititle\fi}}
\immediate\write\gtoutfile{Subj-class: GT or SG, GR etc}
\immediate\write\gtoutfile{MSC-class: \theprimaryclass\ifx\thesecondaryclass\relax\else, \thesecondaryclass\fi}
\immediate\write\gtoutfile{Journal-ref: Algebr. Geom. Topol. \thevolumenumber\s
(\thevolumeyear) \startpage-\finishpage}
\immediate\write\gtoutfile{Comments: Published by Algebraic and
Geometric Topology at}
\immediate\write\gtoutfile{\s\s\s  http://www.maths.warwick.ac.uk/agt/AGTVol\thevolumenumber/agt-\thevolumenumber-\thepapernumber.abs.html}
\immediate\write\gtoutfile{\noexpand\\}
\immediate\write\gtoutfile{}
\ifx\theasciiabstract\relax
\immediate\write\gtoutfile{\theabstract}\else
\immediate\write\gtoutfile{\theasciiabstract}\fi
\immediate\write\gtoutfile{}
\immediate\write\gtoutfile{\noexpand\\}
\immediate\write\gtoutfile{}
\immediate\closeout\gtoutfile}}  

\def\maketitlepage{\makeagttitle\makeheadfile}
\let\makeshorttitle\maketitlepage
\let\maketitle\maketitlepage

\lognumber{36}
\volumenumber{5}
\volumeyear{2005}
\papernumber{36}
\pagenumbers{865}{897}
\received{23 July 2004} 
\accepted{9 May 2005}
\published{29 July 2005}

\usepackage{amsmath,amssymb}
\usepackage{epsfig}
\allowdisplaybreaks

 \newcommand{\Z}{{\mathbb Z}}
 \newcommand{\C}{{\mathbb C}}
 \newcommand{\A}{{\mathbb A}}
 \newcommand{\Q}{{\mathbb Q}}
 \newcommand{\R}{{\mathbb R}}
 \newcommand{\K}{{\mathbb K}}

 \newtheorem{theorem}{Theorem}
 \newtheorem{lemma}[theorem]{Lemma}
  \newtheorem{proposition}[theorem]{Proposition}
 \newtheorem{corollary}[theorem]{Corollary}

 \newtheorem{definition}[theorem]{Definition}
 
 \newtheorem{lemmadef}[theorem]{Lemma + Definition}

\newcommand{\lb}{\left\langle}
\newcommand{\rb}{\right\rangle}
\newcommand{\diag}[2]{\parbox{#2}{\psfig{figure=#1.eps,height=#2}}}
\renewcommand{\o}{{\mathcal{O}}}
\newcommand{\wt}{\widetilde}
\newcommand{\comment}[1]{}

\title{Skein theory for $SU(n)$-quantum invariants}
\asciititle{Skein theory for SU(n)-quantum invariants}

\author{Adam S. Sikora}
\address{Department of Mathematics, University at 
Buffalo\\Buffalo, NY 14260-2900, USA}
\email{asikora@buffalo.edu}

\begin{abstract}
For any $n\geq 2$ we define an isotopy invariant, $\lb \Gamma\rb_n,$ 
for a certain set of $n$-valent ribbon graphs $\Gamma$ in $\R^3,$
including all framed oriented links. We
show that our bracket coincides with the Kauffman bracket for $n=2$
and with the Kuperberg's bracket for
$n=3.$ Furthermore, we prove that for any $n,$ our bracket of a link
$L$ is equal, up to normalization, to the $SU_n$-quantum invariant of $L.$
We show a number of properties of our bracket extending those
of the Kauffman's and Kuperberg's brackets, and we relate it to the
bracket of Murakami-Ohtsuki-Yamada. Finally, on the basis of 
the skein relations satisfied by $\lb\cdot \rb_n,$
we define the $SU_n$-skein module of any $3$-manifold $M$ and we prove
that it determines the $SL_n$-character variety of $\pi_1(M).$
\end{abstract}

\asciiabstract{%
For any n>1 we define an isotopy invariant, <Gamma>_n,$ for a certain
set of n$valent ribbon graphs Gamma in R^3, including all framed
oriented links.  We show that our bracket coincides with the Kauffman
bracket for n=2 and with the Kuperberg's bracket for n=3.
Furthermore, we prove that for any n, our bracket of a link L is
equal, up to normalization, to the SU_n-quantum invariant of L.  We
show a number of properties of our bracket extending those of the
Kauffman's and Kuperberg's brackets, and we relate it to the bracket
of Murakami-Ohtsuki-Yamada. Finally, on the basis of the skein
relations satisfied by <.>_n, we define the SU_n-skein module of any
3-manifold M and we prove that it determines the SL_n-character
variety of pi_1(M).}


\begin{document}

\maketitle

\section{Introduction} 

The $SU_2$-quantum invariant of links, known as the Jones polynomial, can
be conveniently defined in terms the Kauffman bracket invariant, \cite{Ka}.
This approach has several advantages, for example, 
leading to definitions of skein modules and Khovanov
homology\footnote{The Kauffman bracket skein relations allow a
  particularly simple definition of Khovanov's $SU_2$-homology groups, 
\cite{Vi}.} -- two notions in the center of current active research --  
see for example \cite{Bu,FGL,Ga,Ge,GL,PS,S2} and \cite{APS, BN, Go,HK, Ja,
K1,K2,KR,Le,Ra,Vi}.
In \cite{Ku}, Kuperberg constructs a bracket isotopy invariant of links
and $3$-valent graphs in $\R^3,$ with properties analogous to those of the
Kauffman bracket, and shows that it coincides with the
$SU_3$-quantum invariant. We extend his work, by
defining a bracket isotopy invariant $\lb \cdot \rb_n$ for any
$n\geq 2$ and by showing that it determines the $SU_n$-quantum invariant.
More specifically, for any $n\geq 2$ we consider the set ${\cal
W}_n(\R^3)$ of {\em $n$-webs} which are ribbon graphs $\Gamma$ in 
$\R^3$ whose coupons are either $n$-valent sources or $n$-valent sinks.
In particular, ${\cal W}_n(\R^3)$ contains all oriented framed links
in $\R^3$ for any $n.$
We define a bracket isotopy invariant of $n$-webs, $\lb \Gamma \rb_n,$ 
and show that it coincides with the Kauffman bracket for $n=2$,
and with the Kuperberg's bracket for $n=3.$

For reader's convenience, we state three different definitions of
$\lb\cdot \rb_n$: by skein relations, (Theorem \ref{skein_def}), 
by a state sum formula, (Proposition \ref{extra_rel}), and as a 
contraction of tensors, (Section \ref{sqdefinition}). 
Furthermore, we show that for any $n,$ $\lb\Gamma \rb_n$ defines
the $SU_n$-quantum invariant of $\Gamma$ with edges of $\Gamma$ labeled by the
defining $SU_n$-representation and the sinks and the sources of
$\Gamma$ labeled by
the $q$-antisymmetrizer and its dual. The proofs are based on
\cite{RT}.

We prove a number of properties of our bracket which extend 
those of the Kauffman's and Kuperberg's brackets.
In particular, $\lb \cdot \rb_n$ satisfies a skein relation which
relates a crossing to its two ``smoothings'',  
cf.\ Proposition \ref{extra_rel}.
Furthermore, there is a state sum formula for $\lb\cdot \rb_n$, 
Theorem \ref{statesum}, which has the ``positivity''
property analogous to that used in the construction of Khovanov and 
Khovanov-Rozansky homology groups, \cite{K1, K2, KR},  cf.\ Proposition
\ref{non-negative}. Our bracket can be used for an alternative 
definition of Khovanov-Rozansky homology groups; cf.\ 
Section \ref{s_sing_links}.

There exists an alternative generalization of the Kauffman bracket 
due to Murakami, Ohtsuki, and Yamada. Their bracket is 
defined for certain $3$-valent colored graphs with a flow,
\cite{MOY,Mu}. It is expressed in terms of our bracket in Sections
\ref{s_sing_links} and \ref{s_MOY}. We believe that our bracket can be 
expressed in terms of Murakami-Ohtsuki-Yamada bracket as well.
Nonetheless, both approaches have their advantages. Perhaps, an 
advantage of our approach is that $\lb \cdot \rb_n$ 
is related directly to the representation theory of $U_q(sl_n),$ 
and that for $q=1$ our skein relations are equivalent to the 
relations between characters of 
$SL(n)$-representations. Furthermore, our relations seem to be
the most appropriate for the definition of $SU_n$-skein modules of
$3$-manifolds, cf.\ Section \ref{s_skeinmod_def}. (Our definition agrees with
those of Ohtsuki and Yamada, \cite{OY}, and Frohman and Zhong,
\cite{FZ}, for $n=3.$)
Several important properties of the Kauffman bracket skein modules have their
generalizations to the $SU_n$-skein modules for any $n.$
In this paper, we show that $SU_n$-skein module 
of a manifold $M$ for $t=1$ is a commutative ring isomorphic to the 
coordinate ring of the $SL_n$-character variety of $\pi_1(M).$ 
We postpone further study of the $SU_n$-skein modules to a 
forthcoming paper.

{\bf Acknowledgments}\qua The author was partially sponsored by 
NSF grants DMS-0307078 and DMS-0111298.

\subsection{Webs}\label{s_webs} 
An {\it $n$-web} is a ribbon graph in $\R^3,$ cf.\ \cite{RT},
whose every coupon is either an $n$-valent sink or an $n$-valent source.
We denote the coupons of the ribbon graphs by discs rather than
rectangles and we use a marking point to represent the side of 
the coupon with no bands attached,
$\diag{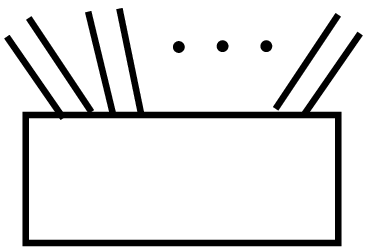}{.3in}\hspace*{.15in}=\diag{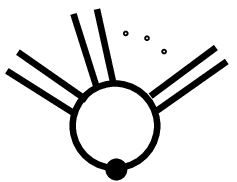}{.3in}\hspace*{.15in}.$

\begin{figure}[ht!]
\centerline{\diag{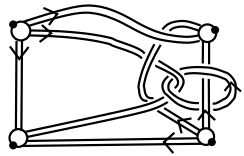}{.8in}\hspace{.3in}}
\caption{An example of a $3$-web}
\end{figure}

For reader's convenience, we restate the definition 
of a web without invoking the notion of a ribbon graph.
The role of ``edges'' of webs is played by {\it bands} which are embeddings
of squares $[0,1]\times [0,1]$ into $\R^3.$
The segments $[0,1]\times \{0\}$ and $[0,1]\times \{1\}$ are
{\it the source} and {\it the target} of the band, respectively. 
Their complement, $[0,1]\times (0,1),$ is {\it the interior} of the band.

An {\it $n$-web} is an oriented surface embedded in $\R^3$ composed of
a finite number of annuli, bands, and discs satisfying
the following conditions:
\begin{enumerate}
\item[(i)] The annuli, disks, and the interiors of the bands are disjoint
from each other.
\item[(ii)] The sources and the targets of bands are disjoint from each other
and all of them lie in the boundaries of discs.
\item[(iii)] The boundary of every disk contains either precisely $n$ sources 
and no targets of bands, in which case the disc is called {\it a source}, or 
it contains precisely $n$ targets and no sources of bands.
In the that case, the disk is called a sink of the web.
\item[(iv)] The marked boundary points of disks lie outside 
the sources and targets of bands.
\end{enumerate}

Since each $n$-web retracts to its spine, which is
an oriented graph, often the bands and discs of webs will be called
its edges and vertices, respectively. In this terminology,
each vertex $v$ of a web is $n$-valent and all edges
adjacent to $v$ are either directed outwards, if $v$ is a source, or 
inwards, if $v$ is a sink. Notice that each web has an equal number of 
sources and sinks.

Our definition of $n$-webs is modeled on the notion of $n$-valent graphs
considered in \cite{S1}, cf.\ Section \ref{scharvar}. The $n$-webs extend the
notion of webs for the geometric $A_1$-spider introduced in \cite{Ku}, 
cf.\ Section \ref{s_kuperberg}.

By analogy with the notion of a link diagram, 
{\it an $n$-web diagram} is a projection $\pi:\Gamma \to \R^2$ of an 
$n$-web $\Gamma$ into $\R^2$ which is an 
embedding of $\Gamma$ except for a finite set of transverse (double)
intersections of bands of $\Gamma$ called crossings.
We require that $\pi$ preserves the orientation of $\Gamma$
(considered as an oriented surface)
and that it embeds the sinks and the sources into $\R^2$ away from the
intersections. In particular, unlike in \cite{RT}, a web diagram is not
allowed to have twists, \diag{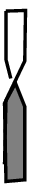}{.3in}\hspace*{-.2in}, in their bands.  
Each web $\Gamma$ is represented by a web diagram; for example:
\centerline{\diag{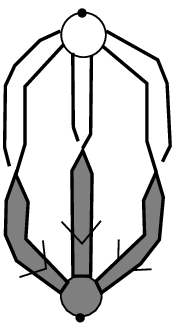}{.8in}\hspace*{-.2in} $=$\quad \diag{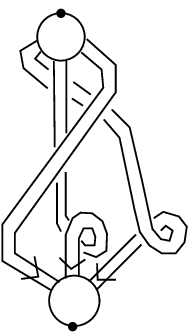}{.8in}\hspace*{-.2in}}

\subsection{The bracket isotopy invariant of $n$-webs}
\label{s_skein_def}
For any permutation $\sigma\in S_n,$ define {\it the length of $\sigma$,} 
$l(\sigma),$ to be the minimal number of factors in the decomposition of 
$\sigma$ into elementary transpositions $(i,i+1),$ $i=1,...,n-1,$
\begin{equation}\label{length}
l(\sigma)=\#\{(i,j): 1\leq i<j\leq n,\ \sigma(i)>\sigma(j)\}.
\end{equation}
For $\sigma\in S_n,$ let \diag{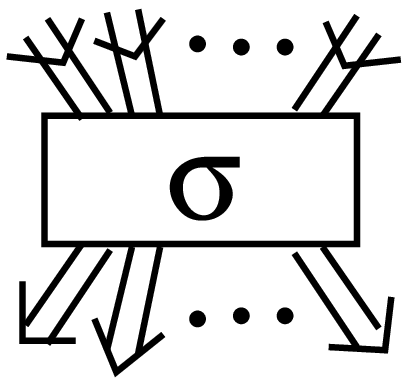}{.5in} denote 
the positive braid with $l(\sigma)$ crossings representing $\sigma.$ 
Such braid is unique. Let $[n]=\frac{q^n-q^{-n}}{q-q^{-1}}$
and let $[n]!=[1]\cdot ...\cdot [n].$

\begin{theorem}\label{skein_def}
There exists a unique isotopy invariant of $n$-webs, $\lb\Gamma \rb_n\in
\Z[q^{\pm \frac{1}{n}}],$  
satisfying the following conditions:
\begin{enumerate}
\item[\rm(i)] $q^\frac{1}{n}\lb \diag{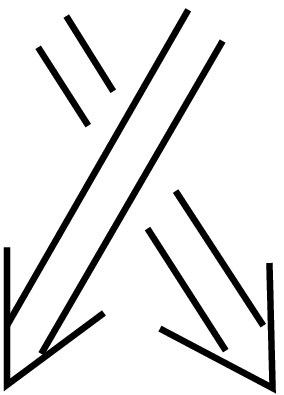}{.3in}\rb_n - q^{-\frac{1}{n}}
\lb\diag{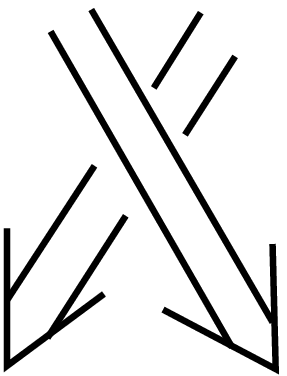}{.3in}\rb_n= (q-q^{-1})\lb\diag{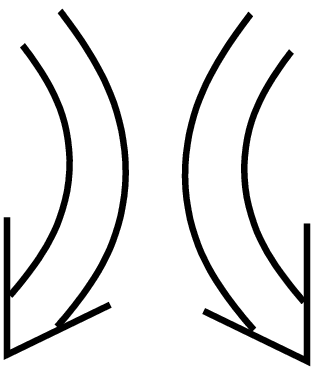}{.3in}\rb_n$
\item[\rm(ii)] $\lb\diag{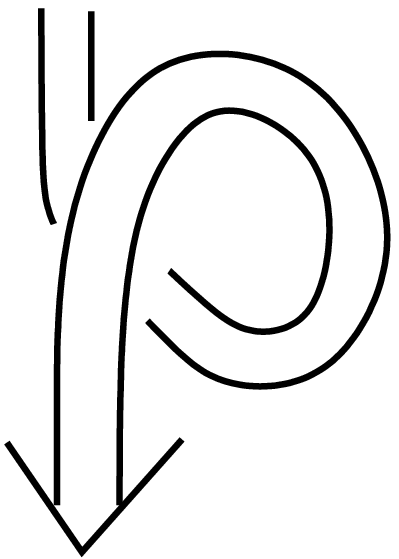}{.3in}\rb_n=q^{n-n^{-1}}\lb
\parbox{.2in}{\psfig{figure=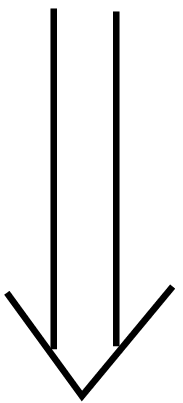,height=.3in}}\rb_n,$
$\lb \diag{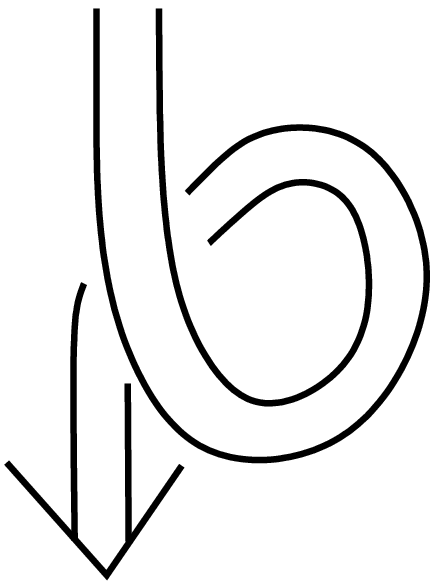}{.3in}\rb_n=q^{n^{-1}-n}\lb 
\parbox{.2in}{\psfig{figure=straight.eps,height=.3in}}\rb_n,$
\item[\rm(iii)] $\lb \diag{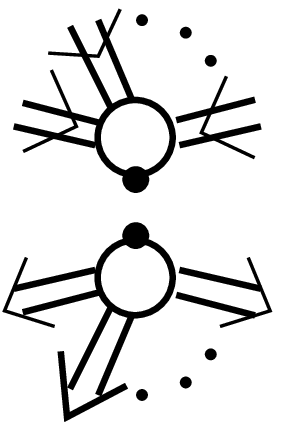}{.5in}\hspace*{-.1in}\rb_n \hspace*{-.1in} =
q^{n(n-1)}\cdot \sum_{\sigma\in S_n} 
(-q^\frac{1-n}{n})^{l(\sigma)}\lb \diag{sigma}{.5in}\ \rb_n.$
\item[\rm(iv)] $\lb \Gamma \cup \diag{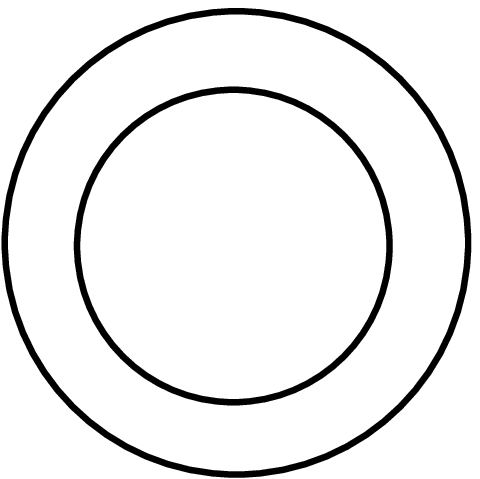}{.2in}\rb_n=[n]\lb\Gamma\rb_n.$
Here $\diag{annulus}{.2in}$ denotes the trivial framed knot 
unlinked with $\Gamma.$
\item[\rm(v)] $\lb \emptyset \rb_n=1$ and, consequently, 
$\lb \diag{annulus}{.2in} \rb_n=[n].$
\end{enumerate}
\end{theorem}

\begin{proof}
The hard part of the statement -- the existence of the bracket --
follows from Theorem \ref{main} stated in Section \ref{sqdefinition}.
The uniqueness of the bracket follows from the fact that each web 
$\Gamma$ has an equal number of sinks and sources, and, therefore, 
condition (iii) makes possible
to represent $\lb \Gamma\rb_n$ by a linear combination of brackets of 
framed links. 
On the other hand, the bracket for framed links is determined 
by conditions (i),(ii), (iv) and (v).
\end{proof}

Relation (iii) appeared in an implicit form in \cite{Bl,Yo} already.

The skein relations of Theorem \ref{skein_def}, appear in the most
natural, but not necessarily, the simplest form.
If $w(\Gamma)$ denotes the writhe (ie. the sum of signs of 
crossings) of a web diagram $\Gamma,$ 
and $v(\Gamma)$ is the number of sinks of $\Gamma$ then
$$P_n(\Gamma)=q^{(n^{-1}-n)w(\Gamma)-n(n-1)v(\Gamma)}\lb \Gamma\rb_n$$ 
is invariant under all Reidemeister moves. Furthermore, it satisfies
the standard skein relations of the $SU(n)$-quantum invariants, cf.\ 
\cite[Thm 4.2.1]{Tu}:

\begin{itemize}
\item $q^nP_n\left(\diag{l+}{.3in}\right)-q^{-n}P_n\left(\diag{l-}{.3in}
\right)=(q-q^{-1}) P_n\left(\diag{l0}{.3in}\right),$
\item $P_n\left(L\cup \diag{annulus}{.2in}\right)=[n]P_n(L).$
\end{itemize}

and the additional relation:
$$\bullet\ P_n\left(\diag{sinksource}{.5in}\hspace*{-.1in}\right) 
\hspace*{-.1in} = \sum_{\sigma\in S_n} 
(-q^{n-1})^{l(\sigma)}P_n\left(\diag{sigma}{.5in}\ \right).$$

\begin{proposition}[Proof in Section \ref{p_skein_def}]\label{extra_rel}
\item $$\lb\diag{l+}{.3in}\rb_n=q^{\frac{n-1}{n}}\lb\diag{l0}{.3in}\rb_n- 
q^{-\frac{n(n-1)}{2}-\frac{1}{n}}\frac{1}{[n-2]!}
\lb \diag{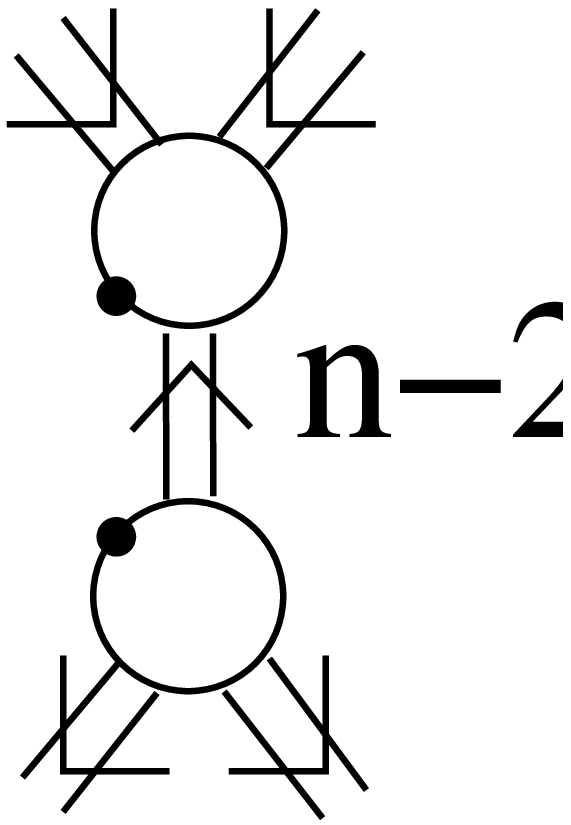}{.5in}\rb_n,$$
where the band labeled by $n-2$ represents $n-2$ parallel bands.
\end{proposition}

The above relation generalizes the Kauffman
bracket skein formula and 
it makes possible to represent any link (or web) as a
linear combination of webs with no crossings. A state-sum formula for
the bracket of webs with no crossings is provided in Section 
\ref{s_state_sum}. Note that various renormalizations of $\lb \cdot
\rb_n$ are possible, leading to a skein formula of Proposition
\ref{extra_rel} without fractional coefficients. Nonetheless, our
definition seems to be the most natural one, cf.\ Section
\ref{sqdefinition}, and leading to the simplest state sum formula.

The following result shows that the bracket 
$\lb \Gamma\rb_n$ for $n$ odd does not depend on the choice of marked 
points on the vertices of $\Gamma.$

\begin{proposition}[Proof in Section \ref{s_marked_proof}]\label{marked}
If $\Gamma,$ $\Gamma'$ are $n$-webs which differ by the choice of
marked points on the boundaries of their discs (vertices) only, then
\begin{enumerate}
\item[\rm(i)] $\lb \Gamma \rb_n=\lb \Gamma' \rb_n$ if $n$ is odd, and
\item[\rm(ii)] $\lb \Gamma \rb_n=\lb \Gamma' \rb_n$ mod $2$ if $n$ is even. 
\end{enumerate}
\end{proposition}

\subsection{The Kauffman bracket and $\lb \cdot \rb_2$}

The Kauffman bracket $[L]\in \Z[A^{\pm 1}]$ is an invariant of unoriented 
framed links $L\subset S^3$ given by the following skein conditions:
$$\left[\diag{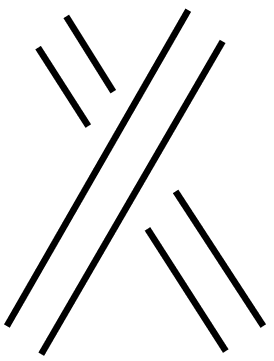}{.25in}\right]=A
\left[\diag{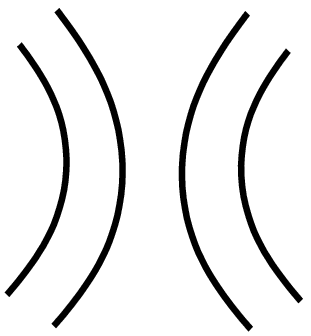}{.2in}\right]+A^{-1} \left [\diag{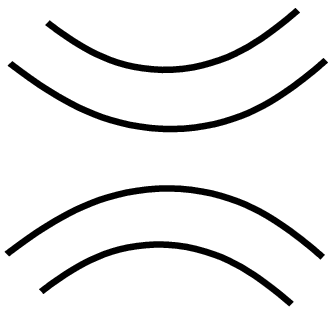}{.2in}\right],
\quad \left[L\cup \diag{annulus}{.2in}\right]=
(-A^2-A^{-2})[L],\quad [\emptyset]=1.$$

\begin{theorem}\label{kb}
For any $2$-web diagram $D,$
$$\lb D \rb_2=(-1)^{w(D)+c(D)}[D],$$
where $A=q^\frac{1}{2},$ $w(D)$ denotes the sum of signs of crossings of 
$D$ and $c(D)$ is the number of components of the link represented by
$D.$ (On the right side $D$ is considered as an unoriented framed 
link diagram).
\end{theorem}

\begin{proof}
The bracket $\lb\cdot \rb_2$ for links is uniquely determined by
conditions (i),(ii),(iv) and (v) of Theorem \ref{skein_def}.
Since $(-1)^{w(D)+c(D)}[D]$ satisfies these relations, the statement follows.
\end{proof}

Note that the bracket of any $2$-web $\Gamma$ can be expressed by 
the bracket of a framed link
by the following operations:
$$\lb \parbox{.3in}{\psfig{figure=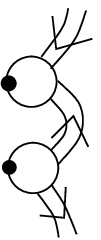,height=.5in}}
\rb_2=q\lb\parbox{.2in}{\psfig{figure=straight.eps,height=.3in}}\rb_2,
\quad \lb \parbox{.3in}{\psfig{figure=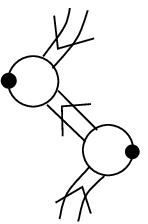,height=.5in}}
\rb_2=-q\lb 
\parbox{.2in}{\psfig{figure=straight.eps,height=.3in}}\rb_2.$$
These equations follow from Theorem \ref{skein_def}.
For example, 
$$\lb \parbox{.3in}{\psfig{figure=kink+-.eps,height=.5in}}
\rb_2=\lb \parbox{.4in}{\psfig{figure=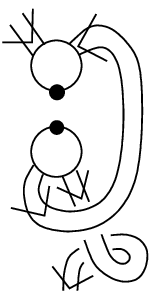,height=.5in}}
\rb_2=q^2\left(
\parbox{.4in}{\psfig{figure=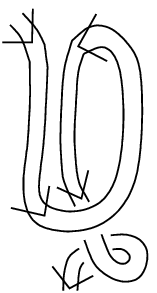,height=.6in}}
-q^{-\frac{1}{2}}
\parbox{.4in}{\psfig{figure=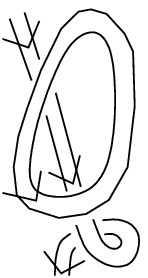,height=.6in}}
\right)=$$
$$q^2\left(q^{-\frac{3}{2}}q^{-\frac{3}{2}}-q^{-\frac{1}{2}}q^{-\frac{3}{2}}
(q+q^{-1})\right)\lb 
\parbox{.2in}{\psfig{figure=straight.eps,height=.3in}} \rb_2=-q \lb 
\parbox{.2in}{\psfig{figure=straight.eps,height=.3in}} \rb_2.$$ 

\subsection{Kuperberg's bracket and $\lb \cdot \rb_3$}
\label{s_kuperberg}

Kuperberg defined an invariant of framed graphs which
are defined as our $3$-webs but without marked points on their vertices,
\cite{Ku}. His bracket is defined by the following relations, in 
which we substituted his $q$ by $q^{-2}$:

\begin{enumerate}
\item[(i)] $\diag{l+}{.3in}= q^{-\frac{1}{3}} \diag{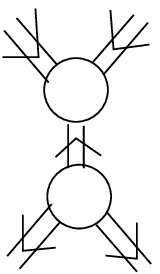}{.4in}+ q^{\frac{2}{3}}
 \diag{l0}{.3in}$
\item[(ii)] $\diag{l-}{.3in}= q^{\frac{1}{3}} \diag{H}{.4in}+ q^{-\frac{2}{3}}
 \diag{l0}{.3in}$
\item[(iii)] $\diag{annulus}{.2in}=[3]$
\item[(iv)] $\diag{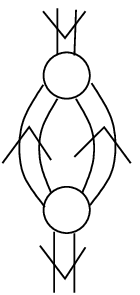}{.4in}=-[2] 
\parbox{.2in}{\psfig{figure=straight.eps,height=.3in}}$
\item[(v)] $\diag{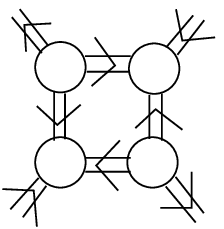}{.4in}=\diag{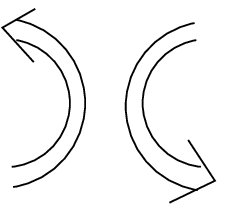}{.3in}+\diag{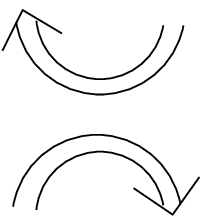}{.3in}.$
\end{enumerate}

Additionally, (i) and (ii) imply\\

\begin{enumerate}
\item[(vi)] $q^\frac{1}{3}\diag{l+}{.3in}-q^{-\frac{1}{3}}\diag{l+}{.3in}=
(q-q^{-1})\diag{l0}{.3in}$
\end{enumerate}
and (i) and (iv) imply
\begin{enumerate}
\item[(vii)] $\diag{kink+}{.3in}=q^\frac{8}{3}
\parbox{.2in}{\psfig{figure=straight.eps,height=.3in}}$.
\end{enumerate}

\begin{theorem}\label{Kuperberg}
Kuperberg's bracket of any web $\Gamma$ is equal to 
$(-q)^{-\frac{3}{2}v(\Gamma)}\lb \Gamma\rb_3,$
where $v(\Gamma)$ is the number of $3$-valent vertices of $\Gamma.$
(By Theorem \ref{marked}(i), $\lb \Gamma\rb_3$ is 
well defined.)
\end{theorem}

\begin{proof}
It is straightforward to check that $q^{-v(\Gamma)}\lb \Gamma\rb_3$
satisfies relations (i),(iii), (vi), and (vii). These equations uniquely
determine Kuperberg's bracket:
(i) makes possible to express Kuperberg's bracket of every
Kuperberg's web as a linear combination Kuperberg's brackets of framed
links. These are uniquely determined by (i),(iii),(vi), and (vii). 
\end{proof}

%
\subsection{Bracket isotopy invariant of framed singular links}
\label{s_sing_links}
%

{\em A singular framed link} is a ribbon graph whose each vertex has two
sinks and two sources. In particular, every oriented framed link
is singular.

\begin{figure}[ht!]
\centerline{\parbox{1in}{\psfig{figure=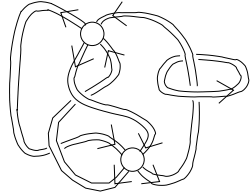,height=1in}}}
\caption{Singular framed link with $2$ singularities}
\end{figure} 

There is a map $$\Psi_n:\{\text{singular framed links in $\R^3$}\}\to
\{\text{$n$-webs in $\R^3$}\},$$
replacing each vertex in a singular framed link by two 
$n$-valent vertices connected by $n-2$ parallel edges:
$$\diag{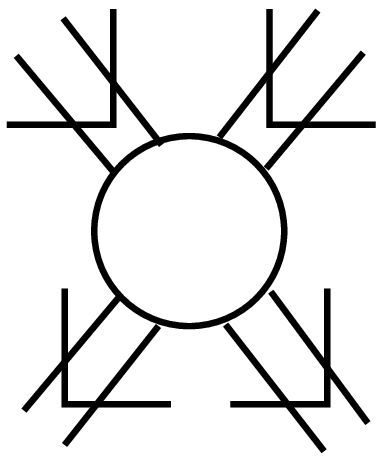}{.3in} \to \diag{crossres}{.5in}$$
For any singular link diagram $D$, let
$$(D)_n= \frac{\lb \Psi(D)\rb_n}{\left(
  [n-2]!q^{n(n-1)/2}\right)^{v(D)}}\cdot q^{w(D)/n},$$
where $v(D)$ is the number of singularities of $D$ (ie. $4$-valent
vertices) and $w(D)$ is the number of positive crossings minus 
negative crossings.

\begin{theorem}\label{sing_links}
$(\Gamma)_n\in \Z[q^{\pm n}],$ and
\begin{enumerate}
\item[\rm(i)] $\left(\diag{l+}{.3in}\right)_n=q\left(\diag{l0}{.3in}\right)_n- 
\left(\diag{fourvertex}{.3in}\right)_n,$
\item[\rm(ii)]
  $\left(\diag{l-}{.3in}\right)_n=q^{-1}\left(\diag{l0}{.3in}\right)_n- 
\left(\diag{fourvertex}{.3in}\right)_n,$
\item[\rm(iii)] $\left(\diag{kink+}{.3in}\right)_n=q^n\left(\diag{straight}{.3in}
\right)_n,$
$\left(\diag{kink-}{.3in}\right)_n=q^{-n}\left(\diag{straight}{.3in}\right)_n,$
\item[\rm(iv)] $\left(L \cup \diag{annulus}{.2in}\right)_n=[n](L)_n.$
\item[\rm(v)] $\left(\emptyset \right)_n=1$ and, consequently, 
$\left( \diag{annulus}{.2in} \right)_n=[n].$
\end{enumerate}
\end{theorem}

\begin{proof} (i) follows from Proposition \ref{extra_rel}.
(ii) follows from (i) and Theorem \ref{skein_def}(i). Parts
(iii)--(v) follow from Theorem \ref{skein_def}(ii),(iv) and (v).
\end{proof}

$(\cdot)$ is a version of the Kauffman-Vogel bracket,
\cite{KV}, Furthermore, it is related to the Murakami-Ohtsuki-Yamada bracket,
\cite[\S 3]{MOY} (see also \cite{Mu}) in the following manner: 
Given a singular framed
link $L,$ label all its edges by $1$ and replace each of its vertices
by a pair of $3$-valent vertices,
$$\diag{fourvertex}{.3in} \to \diag{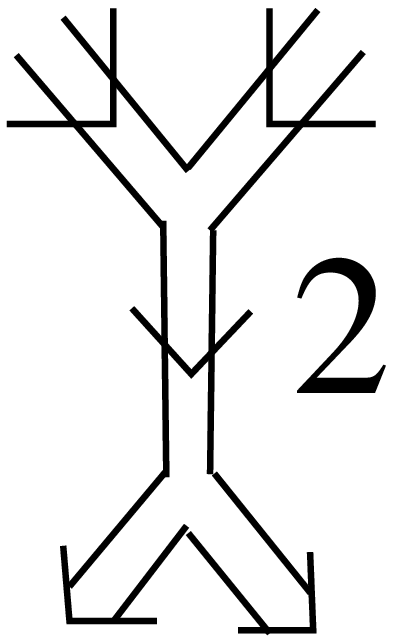}{.5in}.$$
Denote the colored ribbon graph obtained in this way by $\Phi(L).$

\begin{proposition}
$(L)_n$ is equal to the Murakami-Ohtsuki-Yamada bracket of $\Phi(L)$
when our $q$ is identified with $q^\frac{1}{2}$ in \cite{MOY}.
\end{proposition}

\begin{proof} It follows from \cite{MOY} that the 
MOY bracket of $\Phi(L)$ satisfies conditions (i)-(v) of Theorem 
\ref{sing_links}. These conditions determine $(\cdot)$ uniquely.
\end{proof}

The above proposition relates the Murakami-Ohtsuki-Yamada bracket with
our bracket for some graphs only. We will see in Section \ref{s_MOY}, 
that Murakami-Ohtsuki-Yamada bracket of every $3$-valent graph with a 
flow is determined by our bracket of a corresponding $n$-web.

Khovanov and Rozansky use $(\cdot)$ to define a
homology theory whose extended Euler characteristic is the
$SU_n$-quantum invariant, \cite{KR}.

%
\subsection{State sum formula for the brackets of planar webs}
\label{s_state_sum}
%

An important future of Kauffman bracket is that it is given 
by a simple state sum formula. We describe a generalization of this
formula for our bracket of $n$-webs below.
An $n$-web diagram $\Gamma$ is planar if it has no crossings.
Since Proposition \ref{extra_rel} makes possible to express the
bracket of any $n$-web as a linear combination of brackets of planar 
$n$-webs, we formulate a state sum formula planar webs only.

{\it A state} $S$ of a planar $n$-web diagram $\Gamma$ is a labeling of 
its annuli and bands $e$ by numbers $S(e)\in \{1,...,n\}$ 
such that the bands attached to every disc are labeled by 
different numbers. (There is no restriction on labeling of annuli.)

Note that every state of $\Gamma$ determines an ordering of
edges adjacent to every vertex $v$ of $\Gamma.$
However, there is also a natural ordering of edges adjacent to $v,$
which does not depend on the choice of a state:
If $v$ is a sink then we order the edges from $1$
to $n$ by starting at the base point of the disc $v$ and then by
moving clockwise around its boundary. If $v$ is a source then we start 
at the base point 
of the disc $v$ and move counter-clockwise around its boundary.\vspace*{.2in}

\begin{figure}[ht!]
\centerline{\diag{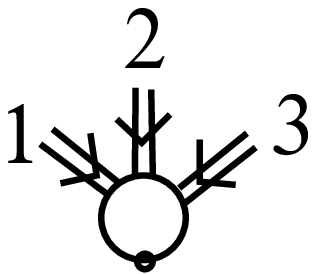}{.4in}\hspace*{.5in} \diag{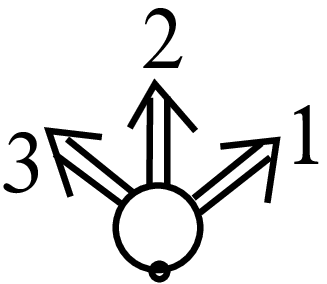}{.4in}}
\caption{The canonical ordering of bands adjacent to
a sink and a source for $n=3$}
\end{figure}

For any state $S$ and a vertex $v$ of $\Gamma,$ let $P(S,v)(i)$
denote the label associated by the state $S$ with the $i$th band 
attached to $v.$ Hence $P(S,v)\in S_n.$

For any state $S,$ we define {\em the rotation index} of $S$ as follows.
\begin{equation}\label{rot}
rot_n(S)=\sum_e ind(a)(2S(e)-n-1),
\end{equation}
where the sum is over all annuli and bands of $\Gamma.$ If $e$ is an 
annulus, then $ind(e)$ is either $+1$ or $-1$ depending on whether 
$e$ is oriented anti-clockwise or clockwise.
The indices $ind(e)$ for edges $e$ of $\Gamma$ are defined as follows:
For each band $e$ in an $n$-web $\Gamma$ choose a smooth embedded arc 
$$\alpha_e:[0,1]\to \text{the band $e$} \cup \text{the sink of $e$}
\cup \text{the source of $e$}$$
connecting the marked points of the sink and 
the source.

Choose the arcs $\alpha_e$ such that for different bands $e,e'$ the 
arcs $\alpha_e,$ $\alpha_{e'}$ are disjoint, except possibly meeting
at one or two of their endpoints. 
The union $\bigcup_e \alpha_e$ taken over all bands $e$ of $\Gamma$
forms an oriented $n$-valent graph $\Gamma'$ in $\R^2.$ 
We say that $\Gamma'$ is {\it a core} of $\Gamma$ if for every 
vertex $v$ of $\Gamma'$ the tangent vectors at $v$ of arcs 
having one of their endpoints at $v$ are pointing in the same
direction.

\begin{figure}[ht!]
\centerline{\diag{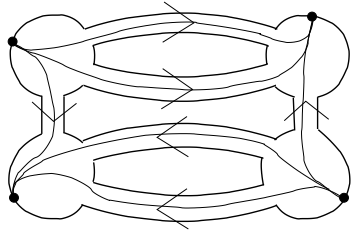}{1.2in}\hspace{.55in}}
\caption{A core in a $3$-web}
\end{figure}

Given a core $\Gamma'$ of $\Gamma,$ for every band $b$ in $\Gamma$
we define its winding number,
$$ind(b)=\frac{1}{2\pi i}\int_0^1 
\frac{\alpha_b''(t)}{\alpha_b'(t)}dt,$$
where we identify $\R^2$ with $\C$ and assume that 
$\alpha_b:[0,1]\to \C.$
Note that $ind(b)=\frac{\beta}{2\pi}$ mod $\Z,$ where $\beta$ is the angle 
between the tangent vectors to $\alpha_b$ at its endpoints.

\begin{lemma}[Proof in Section \ref{p_statesum}]\label{core_indep}
For any state $S$ of an $n$-web $\Gamma,$ $rot_n(S)$ is independent of
the choice of a core of $\Gamma.$ Furthermore, $rot_n(S)$ is an
isotopy invariant of $\Gamma$ and $rot_n(S)\in \Z.$
\end{lemma}

\begin{theorem}[Proof in Section \ref{p_statesum}]\label{statesum}
For any planar $n$-web $\Gamma,$
$$\lb \Gamma \rb_n=\sum_S q^{rot_n(S)} \prod_ v 
(-q)^{l(P(S,v))},$$
where the sum is taken over all states of $\Gamma$ and the product is 
over all its vertices.
\end{theorem}

We leave the proof of the following proposition to the reader.

\begin{proposition}\label{non-negative}
Let $\Gamma$ be a planar $n$-web obtained by resolving all crossings
of a link by the skein relation of Proposition \ref{extra_rel},
(In other words, let $\Gamma$ be one of the leaves of the skein tree of $L$).
Then $\sum_v l(P(S,v))$ is even for any state $S$. Consequently, all
coefficients of $\lb \Gamma\rb_n\in \Z[q^{\pm 1}]$ are non-negative.
\end{proposition}

\subsection{Murakami-Ohtsuki-Yamada colored $3$-valent graphs}\label{s_MOY}

Murakami, Ohtsuki, and Yamada defined an $SU_n$-quantum invariant of links
by using $3$-valent graphs with a flow.
Inspired by their work, we say that a $3$-valent oriented, framed graph
embedded into $\R^3$ is a
Murakami-Ohtsuki-Yamada graph (MOY-graph, for short) if 
the edges of $\Gamma$ are labeled by positive integers forming a flow 
on $\Gamma:$ 
$$\diag{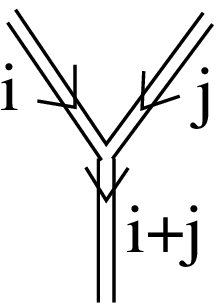}{.6in} \quad \text{or} \quad \diag{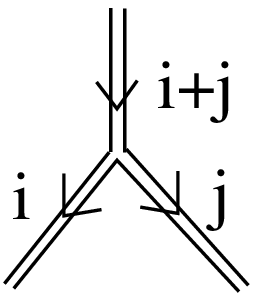}{.6in}$$
We allow annuli embedded into $\R^3$ colored by
positive integers as components of MOY-graphs as well.

An MOY-graph $\Gamma$ is an $MOY_n$-graph if the labels of its
edges and annuli do not exceed $n.$ The purpose of this section is to 
show that our bracket of $n$-webs defines a bracket invariant of 
$MOY_n$-graphs, which coincides (up to a normalization) with the 
Murakami-Ohtsuki-Yamada bracket.

For any $MOY_n$-graph diagram $\Gamma$ 
with no crossings, let $W(\Gamma)$ be a ribbon graph obtained 
by replacing all vertices of $\Gamma$ as follows:
$$\diag{moyfin}{.65in} \to \diag{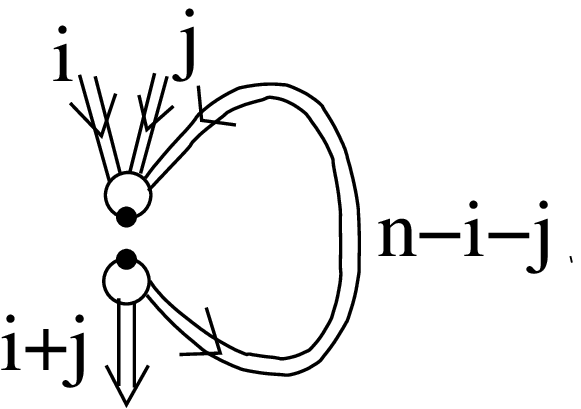}{.7in}\hspace{.3in} ,\qquad
\diag{moyfout}{.65in}\to \diag{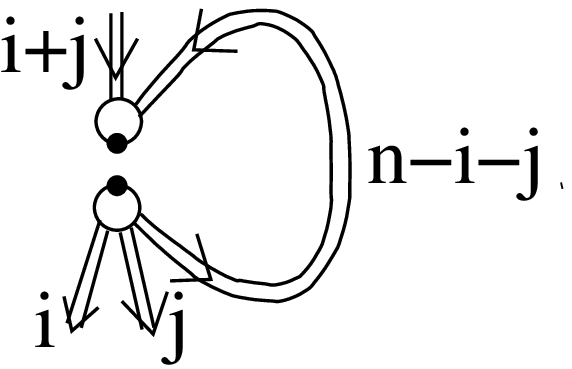}{.65in}$$
As before, an edge of a web labeled by $i$ denotes
$i$ parallel edges. Let
$$[\Gamma]_n=\lb W(\Gamma) \rb_n.$$

\begin{corollary}
$[\Gamma]_n$ is an isotopy invariant of $MOY_n$-graphs.
\end{corollary}

We are going to show that $[\Gamma]_n$ is a renormalization of 
the Murakami-Ohtsuki-Yamada bracket of $\Gamma.$

For any MOY-graph $\Gamma$ denote the labels of edges $e$ of $\Gamma$ by $|e|.$
For any vertex $v$ of $\Gamma,$ denote the adjacent edge with the largest
label by $e^v_0,$ and the left and the right of the two other
adjacent edges by $e^v_1$ and by $e^v_2$ respectively. Hence the
adjacent edges to $v$ in $\Gamma$ are either

\begin{center}
\diag{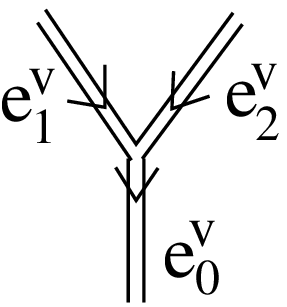}{.6in}\quad or \quad\diag{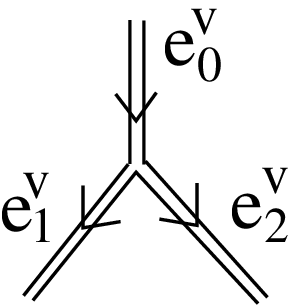}{.6in}
\end{center}

\noindent By the definition of a flow, $|e^v_0|=|e^v_1|+|e^v_2|$ for
any vertex $v.$

We say that a function $s$ assigning an $|e|$-element subset of
$\{1,...,n\}$ to every edge $e$ of $\Gamma$ is {\em an $n$-state} (or,
simply, {\em a state}) of $\Gamma$ if $s(e^v_1)\cap s(e^v_2)=\emptyset$ and 
$s(e^v_1)\cup s(e^v_2)=s(e^v_0)$ for every vertex $v.$

Note that our definition coincides with the definition of \cite{MOY} 
if the sets $s(e)=\{i_1,...,i_{|e|}\}$ and
$\{i_1-\frac{n-1}{2}, i_2-\frac{n-1}{2},
..., i_{|e|}-\frac{n-1}{2}\}$ are identified.

Any $n$-state $s$ splits $\Gamma$ into several
simple closed loops (which may intersect each other), each labeled by
an integer between $1$ and $n.$
Following \cite{MOY}, let the rotation number of an $n$-state $s$ be
$$rot(s)=\sum_C (s(C)-\frac{n+1}{2})rot(C)\in \frac{1}{2}\Z,$$ where 
the sum is over all simple closed loops $C$ of the splitting of
$\Gamma$ by $s,$ $s(C)$ is the label of $C,$ and
$rot(C)$ is either $+1$ or $-1$ depending on whether $C$ is oriented
anti-clockwise or clockwise.

As in \cite{MOY}, for any two sets $s_1,s_2\subset \{1,...,n\}$ 
we denote by $\pi(s_1,s_2)$ the number of pairs 
$(i_1,i_2)\subset s_1\times s_2$ such that $i_1>i_2.$

\begin{proposition}[Proof in Section \ref{sMOY-state_sum_proof}]\label{MOY-state_sum} 
For any $MOY_n$-graph diagram $\Gamma$ with no crossings
$$[\Gamma]_n= {\cal N}(\Gamma)\sum_{\text{$n$-states}\ s} q^{2rot(s)}
\prod_{\text{vertices}\ v} (-q)^{\pi(s(e_1^v),s(e_2^v))}.$$
${\cal N}(\Gamma)$ is a normalization factor,
$${\cal N}(\Gamma)=
\prod_e [|e|]!\cdot
\prod_v q^\frac{n(n-1)-|e^v_1|\cdot|e^v_2|}{2}
[n-|e_0^v|]!$$
where the first product is taken over all edges of $\Gamma$ 
and the second product is over all vertices of $\Gamma.$ (Annuli of
$\Gamma$ are not considered as edges.)
\end{proposition}

In order to avoid confusion with our bracket, we denote the
$n$-th Murakami-Ohtsuki-Yamada bracket of an $MOY_n$-graph $\Gamma$ by 
$\{\Gamma\}_n.$

\begin{proposition}[Proof in Section \ref{sMOY_bracket_proof}]
\label{MOY-bracket}
For any $MOY_n$-graph diagram $\Gamma$ with no crossings,
$$\{\Gamma\}_n=q^{-\frac{1}{4}\sum_v |e_1^v|\cdot |e_2^v|}
\sum_{\text{states}\ s} q^{rot(s)}
\prod_{\text{vertices}\ v} q^{\pi(s(e_1^v),s(e_2^v))/2}.$$
\end{proposition}

The following lemma is needed to relate the brackets
$\{\cdot\}_n$ and $\lb \cdot\rb_n:$

\begin{lemmadef}
For any MOY-graph $\Gamma,$
$$\eta_n(\Gamma)=2rot(s) \mod 2$$ does not depend on the $n$-state $s$ of
$\Gamma.$
\end{lemmadef}

\begin{proof}
Since $2rot(s)=(n+1)\sum_C rot(C)$ mod $2,$
it is enough to show that for any state $s$, the induced splitting of 
$\Gamma$ into simple loops $\{C\}$ is such that $\sum_C rot(C)$ does 
not depend on $s.$ To prove that, consider all cups and caps, $c,$ of 
$\Gamma.$ 
Each of them, being a part of an edge of $\Gamma$, has an associated flow
$|c|.$ Let $rot(c)$ be either $+1/2$ or $-1/2$ depending on whether $c$ is
oriented anti-clockwise or clockwise,
$$\begin{array}{cccc}
\diag{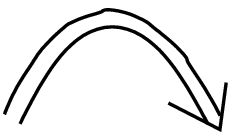}{.2in} & \diag{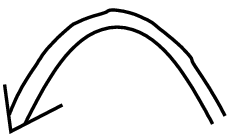}{.2in} & \diag{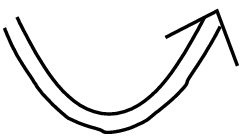}{.2in} &
\diag{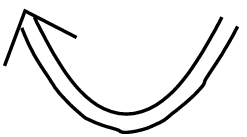}{.2in}\\
rot(c)=-1/2 & rot(c)=1/2 &rot(c)=1/2 &rot(c)=-1/2\\
\end{array}$$
Note that
$$\sum_C rot(C)=\sum_{\text{cups and caps}:\ c} |c|rot(c),$$
and hence the left hand side does not depend on $s.$
\end{proof}

By Proposition \ref{MOY-bracket},
$$\{\Gamma\}_n q^{\frac{1}{4}\sum_v |e_1^v|\cdot |e_2^v|}\in \Z[q^{\pm
    \frac{1}{2}}]$$ and,
by Proposition \ref{MOY-state_sum}, substitution $q^\frac{1}{2}\to -q$ 
gives
$$\left(\{\Gamma\}_n q^{\frac{1}{4}\sum_v |e_1^v|\cdot |e_2^v|}
\right)_{q^\frac{1}{2}\to -q}=\sum_{\text{states}\ s} q^{2rot(s)}
\prod_{\text{vertices}\ v} (-q)^{\pi(s(e_1^v),s(e_2^v))}=$$
$$=[\Gamma]_n\cdot (-1)^{\eta_n(\Gamma)}/{\cal N}(\Gamma).$$

\begin{corollary} The value of the Murakami-Ohtsuki-Yamada bracket of
any $MOY_n$-graph $\Gamma$ is determined by $\lb W(\Gamma)\rb_n.$
\end{corollary}

%
\section{Definition of the bracket using tensors}
\label{sqdefinition}
%

We will state now another, more explicit definition of the bracket
of $n$-webs, which utilizes the construction of Reshetikhin and Turaev, 
\cite{RT}. Given a ribbon Hopf algebra $H,$ they constructed an isotopy
invariant for ribbon graphs whose edges are labeled by representations
of $H$ and whose vertices are labeled by $H$-invariant tensors.
We are going to see that our bracket $\lb \Gamma\rb_n$ 
is the Reshetikhin-Turaev quantum $sl(n)$ invariant for $\Gamma$ considered
as a ribbon graph whose edges are decorated by the defining representation
$V$ and whose sinks and sources and decorated by an
element of the $1$-dimensional representation $\bigwedge^n V\subset
V^{\otimes n}$ and by its dual, respectively.

Let $V$ be an $n$-dimensional vector space over $\C(q)$ with a basis
$e_1,...,e_n.$
Given a web diagram $\Gamma$ decompose it into pieces with the following
tensors associated with them:
\begin{equation}\label{crossings}
\begin{array}{cccc}
\diag{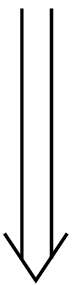}{.3in} & \diag{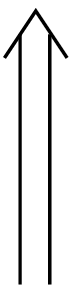}{.3in} & \diag{l+}{.3in} & 
\diag{l-}{.3in}\\
Id_V\ & Id_{V^*}\ & \hat R\ & {\hat R}^{-1},\ \\
\end{array}
\end{equation}
where $\hat R: V\otimes V\to  V\otimes V$ is given by
\begin{equation}\label{R-matrix}
{\hat R}(e_i\otimes e_j)=q^{-\frac{1}{n}}
\begin{cases}
e_j\otimes e_i & \text{if $i>j$,}\\
q e_i\otimes e_j  & \text{if $i=j$,}\\
e_j\otimes e_i +(q-q^{-1})e_i\otimes e_j & \text{if $i<j$.}\\
\end{cases}
\end{equation}
\begin{equation}\label{cupscaps}
\begin{array}{cccc}
\diag{cap+}{.2in} & \diag{cap-}{.2in} & \diag{cup+}{.2in} & \diag{cup-}{.2in}\\
\text{\small $\sum_i e^i\otimes e_j\to \delta_{ij}$} & 
\text{\small $e_i\otimes e^j\to q^{2i-n-1}\delta_{ij}$}&
\text{\small $1\to \sum_i e_i\otimes e^i$} &
\text{\small $1\to \sum_i  q^{n+1-2i} e^i\otimes e_i$}\\
\text{\small $V^*\otimes V\to \C(q)$}& 
\text{\small $V\otimes V^*\to \C(q)$} &
\text{\small $\C(q)\to V\otimes V^*$} & 
\text{\small $\C(q)\to V^*\otimes V$}\\
\end{array}
\end{equation}
\begin{equation}\label{sourcesink}
\begin{array}{cc}
\diag{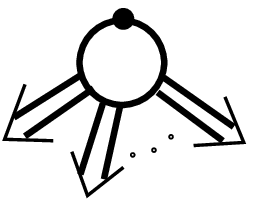}{.3in} & \diag{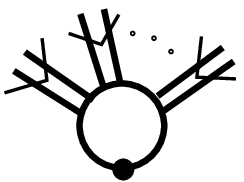}{.3in}\\
\text{\small ${\cal T}_-:V^{\otimes n}\to \C(q)$} & 
\text{\small ${\cal T}_+:\C(q)\to V^{\otimes n}$,}\\
\end{array}
\end{equation}
where 
\begin{equation}\label{T+}
{\cal T}_+(1)=T_+=\sum_{\sigma\in S_n} (-q)^{l(\sigma)} 
e_{\sigma(1)}\otimes e_{\sigma(2)}\otimes ... \otimes e_{\sigma(n)},
\end{equation}
and
\begin{equation}\label{T-}
{\cal T}_-(e_{i_1}\otimes ...\otimes e_{i_n})=
\begin{cases}
(-q)^{l(\sigma)} & \text{if $(1,...,n)\to (i_1,...,i_n)$ is a permutation
$\sigma$}\\
0 & \text{otherwise.}\\
\end{cases}
\end{equation}

\begin{definition}\label{qdef}
For any $n$-web diagram $\Gamma$ decomposed into pieces as above,
let $\lb \Gamma\rb_n\in \Q(q^{\frac{1}{n}})$ be the scalar obtained by
the contraction of the corresponding tensors.
\end{definition}

\begin{theorem}[Proof in Section \ref{s_proof_main}]\label{main}
\begin{itemize}
\item[\rm(i)] The bracket defined above is an isotopy invariant of $n$-webs.\\
\item[\rm(ii)] It satisfies the properties of the bracket stated in 
Theorem \ref{skein_def}.
\end{itemize}
\end{theorem}

Consequently, Definition \ref{qdef} coincides with the definition of 
the bracket given in Theorem \ref{skein_def}.

\section{The $SU(n)$-skein modules of $3$-manifolds}
\subsection{The definition of the skein module}\label{s_skeinmod_def}

Let $M$ be an orientable $3$-manifold, possibly with non-empty
boundary, and let $n\geq 2.$ Let ${\cal W}_n(M)$ denote the set of all 
isotopy classes of $n$-webs embedded into $M,$ including the empty
web, $\emptyset.$
Consider a ring $R$ with a specified invertible 
element $t.$ The $SU(n)$-skein module of $M$ with coefficients in $R$ is
the quotient of the free $R$-module $R{\cal W}_n(M)$
by relations 
\begin{enumerate}
\item[(i)] $t \diag{l+}{.3in} - t^{-1}\diag{l-}{.3in}-
  (t^n-t^{-n})\diag{l0}{.3in}$
\item[(ii)] $\diag{kink+}{.3in}-t^{n^2-1}
\parbox{.2in}{\psfig{figure=straight.eps,height=.3in}},$
$\diag{kink-}{.3in}-t^{1-n^2}
\parbox{.2in}{\psfig{figure=straight.eps,height=.3in}},$
\item[(iii)] $\diag{sinksource}{.5in}\hspace*{-.1in} \hspace*{-.1in}-
t^{n^3-n^2}\cdot \sum_{\sigma\in S_n} 
(-t^{(1-n)})^{l(\sigma)} \diag{sigma}{.5in}.$
\item[(iv)] $\Gamma\cup \diag{annulus}{.2in} -[n]\Gamma.$
\end{enumerate}

Note that the above relations correspond to equations (i)-(v) of 
Theorem \ref{skein_def}, after the substitution $t=q^\frac{1}{n}.$
We denote the quotient module by ${\cal S}_n(M;R,t).$

The $SU_3$-skein module was defined independently in \cite{FZ}, and
earlier, for $3$-dimensional thickenings of surfaces, in \cite{OY}. 
The definitions of Frohman-Zhong and Ohtsuki-Yamada are equivalent to
ours, since their skein relations are the skein relations for
$q^{-v(\Gamma)}\lb \Gamma\rb_n,$ when $A=-q^\frac{1}{3}=-t,$ 
cf.\ Theorem \ref{Kuperberg}.

Theorem \ref{kb} implies the following:

\begin{corollary}
${\cal S}_2(M;R,t)$ is isomorphic to the Kauffman bracket skein module 
of $M$ with coefficients in $R$ and $A=t.$
\end{corollary}

It follows directly from the definition, that if $f:R\to R'$ is a
homomorphism of rings such that $f(t)=t'$ then
$${\cal S}_n(M;R',t')={\cal S}_n(M;R,t)\otimes_R R'.$$
Since for any ring $R$ with an invertible element $t$ there is a 
map $f:\Z[t^{\pm 1}]\to R,$ Theorem \ref{skein_def} can be restated
as follows:

\begin{corollary} For any ring $R$ with an invertible element $t,$
${\cal S}_n(\R^3;R,t)=R$ 
\end{corollary}

Below, we describe the relation between ${\cal S}_n(M;R,t)$ and
$SL_n(R)$-representa\-tions of $\pi_1(M),$ which generalizes the
theorems of Bullock, \cite{Bu}, and ours with J. Przytycki, \cite{PS},
for the Kauffman bracket skein modules.
Further analysis of $SU(n)$-skein modules is postponed to \cite{S2}. 
The discussion below and the results of \cite{S2} show that 
$SU(n)$-skein modules have many properties analogous to those of
of the Kauffman bracket skein module.

%
\subsection{$SU(n)$-skein modules and character varieties}
\label{scharvar}
%

Since the skein relation (i) above reduces to 
$\diag{l+}{.3in} = \diag{l-}{.3in}$ for $t=1,$ the $n$-webs in $M$,
which are freely homotopic to each other, are identified in 
${\cal S}_n(M;R,1).$ Consequently, the operation of taking the
disjoint union, $\Gamma_1,\Gamma_2\to \Gamma_1\cup\Gamma_2,$ 
extends to a well defined product in ${\cal S}_n(M;R,1)$ making this
module a commutative $R$-algebra.

Furthermore, note that as an $R$-algebra,
${\cal S}_n(M;R,1)$ is isomorphic to $\A_n(M)$ (over the ring of 
coefficients $R$) defined in \cite{S1}. Hence, by 
\cite[Theorem 3.7]{S1} we have:

\begin{corollary}
If $\K$ is an algebraically closed field of characteristic $0$ then
${\cal S}_n(M;\K,1)\simeq \o(X_n(\pi_1(M))),$ where 
$\o(X(M))$ denotes the ring of global sections (the coordinate ring)
of the $SL_n(\K)$-character variety of $\pi_1(M).$
\end{corollary}

The $SL_n(\K)$-character variety, $X_n(G),$ of a group $G$ is an affine
algebraic scheme over $\K$ ``describing'' the $SL_n(\K)$-representations of
$G$ up to conjugation.
More precisely, the closed points of $X_n(G)$ (ie. the maximal ideals
in $\o(X(G))$ correspond to the semi-simple $SL_n(\K)$-representations
of $G$ up to conjugation. For a precise definition of $SL_n$-character
varieties see \cite{S1,LM}.

Up to nilpotent elements, the ring $\o(X(G))$ can be described as
follows:
A characteristic function $f:SL_n(\K)\to \K$ is any polynomial in the
entries of the matrices in $SL_n(\K)$ which is invariant under the
conjugation by $SL_n(\K).$ Characteristic functions of $SL_n(\K)$ form
a $\K$-algebra generated by the functions $f_n(A)=tr(A^n).$
Let $X'(G)$ be the set of all {\em generalized $SL_n$-characters} of
$G,$ that is $\K$-valued functions on $G$ of the form $\psi=f\circ\rho,$
where $\rho:G\to SL_n(\K)$ is a representation and $f$ is a
characteristic function on $SL_n(\K).$ With any $g\in G$ there
is the associated ``evaluation at $g$'' function $\tau_g: X'(G)\to \K,$ 
$\tau_g(\psi)=\psi(g).$ The $\K$-algebra $\o(X(G))/\sqrt{0}$ is isomorphic to
the $\K$-algebra generated by all $\tau_g$ for all $g.$ 
If we think of the functions $\tau_g$ as regular functions on
$X_n'(G),$ then for any finitely generated group $G,$ $X_n'(G)$ becomes
an affine algebraic set whose coordinate ring is isomorphic to
$\o(X(G))/\sqrt{0}.$

%
\section{Preliminaries for the proofs}
%


\subsection{The quantum $sl_n$ group and its
defining representation}

Let $U_h=U_h(sl_n)$ be defined as in \cite[Sect. 6.1.3]{KS}.
Note that Reshetikhin's and Turaev's definition of 
$U_h(sl_n),$ \cite[Sect. 7.1]{RT},
coincides with our definition of $U_h$ after taking into account 
the following changes:
\begin{itemize}
\item $X^+_i$ and $X^-_i$ in \cite{RT} are $E_i$ and $F_i$ in \cite{KS},
\item $h$ used by Reshetikhin and Turaev is equal to $2h$ in \cite{KS}.
\end{itemize}
In this paper we will use Klimyk-Schm\"udgen $h$ only.
Let $V=\C[[h]]^n$ be the defining representation of $U_h,$
presented explicitly in \cite[Sect. 8.4.1]{KS}.
Additionally, consider the quantum group
$U_q=U_q(sl_n),$ \cite[Sect. 6.1.2]{KS}, and 
its defining representation, $V=C(q)^n,$ as defined in 
\cite[Sect. 8.4.1]{KS}. This double meaning of $V$ will not lead to 
confusion since the $U_q$ and $U_h$ actions on $V$ agree if
$q$ and $K_i$ are identified with $e^h$ and $e^{hH_i}$ respectively.
In both cases, $e_1,...,e_n$ will be the weight basis of $V$ with 
the heighest weight vector $e_1$, see \cite[Sect. 8.4.1]{KS}.

The defining representation $\rho: U_q\to End(V)$ is given by the matrices
\begin{equation}\label{K-action}
\rho(K_i)=q^{-1}E_{ii}+qE_{i+1,i+1}+\sum_{j\ne i,i+1}E_{jj},\quad 
i=1,2,...,n-1,
\end{equation}
\begin{equation}\label{EF-action}
\rho(E_i)=E_{i+1,i},\quad \rho(F_i)=E_{i,i+1},\quad i=1,2,...,n-1,
\end{equation}
where $E_{ij}$ denotes the matrix whose all entries are $0$ except
for the $(i,j)$th entry which is $1,$ cf.\ \cite[Sect. 8.4.1]{KS}.

\subsection{The Hecke algebra associated with $U_q$}\label{s_hecke}
For the purpose of studying $U_q(sl_n)$-actions on $V^{\otimes k}$
we define the Hecke algebra, $H_k,$ as follows. Let $H_k$ be 
the non-commutative, associative algebra over $\C(t)$
generated by elements $g_1^{\pm 1},...,g_k^{\pm 1}$ subject to the 
following relations:
$g_ig_j=g_jg_i$ for $|i-j|\geq 2,$ $g_ig_{i+1}g_i=g_{i+1}g_ig_{i+1},$ for
$i=1,...,k-1,$ and $(g_i-t^n)(g_i+t^{-n})=0$ for all $i.$ The algebra $H_k$ 
is an $k!$-dimensional space with basis vectors $h_\sigma,$ for 
$\sigma\in S_k,$ which satisfy the following multiplication rules
\begin{enumerate}
\item[(i)] $(h_{(i,i+1)}-t^{n-1})(h_{(i,i+1)}+t^{-n-1})=0$ 
for $i=1,...,n-1,$
\item[(ii)] $h_{\sigma}h_{\tau}=h_{\sigma\tau}$ if $l(\sigma\tau)=
l(\sigma)+l(\tau).$
\end{enumerate}
We have $g_i=t\cdot h_{(i,i+1)}.$
Although various definitions of Hecke algebra appear in the literature, 
see for example \cite[Sect. 8.6.4]{KS}, \cite[Sect. 12.3]{CP}, \cite{Gy}, 
they are all isomorphic to our $H_k$
after a proper extension of base field.
In fact, substituting the quadratic equation (i) above by any other quadratic 
equation with distinct roots leads an isomorphic algebra, after a 
finite extension of the base field. For example, Gyoja's $q$ in 
\cite{Gy} is our $t^{2n}$
and Gyoja's $h(\sigma)$ is ours $t^{(n+1)l(\sigma)}h_\sigma.$
Our somewhat cumbersome notation is chosen so that 
$h_{(i,i+1)}$'s satisfy the same quadratic equation as $\hat R$
for $q=t^n.$ Therefore, from now on, we will assume that $q=t^n.$
We summarize the basic relations between $H_n$ and the defining
representation $V$ of $U_q:$

\begin{proposition}\label{prop_hecke}$\phantom{99}$
\begin{itemize} 
\item[\rm(i)] $H_k$ acts on $V^{\otimes k}$ in such a way
that the actions of $h_{\sigma}$ and of \diag{sigma}{.4in} 
on $V^{\otimes k}$ coincide.
\item[\rm(ii)] The $H_k$ and $U_q$ actions on $V^{\otimes k}$ commute.
\item[\rm(iii)] {\rm(Frobenius-Schur duality)}\qua The images of the maps
$U_q\to End(V^{\otimes k})$ and $H_n\to End(V^{\otimes k})$ are 
centralizers of each other.
In particular, any $U_q$-equivariant endomorphism of $V^{\otimes k}$
is of the form $w\to x\cdot w$ for a certain $x\in H_k.$
\end{itemize}
\end{proposition}

In \cite[page 843]{Gy}, Gyoja defines two elements $e_+,e_-$ which in
our notation are:
\begin{equation}\label{e-def}
e_+=\sum_{\sigma\in S_k} q^{\frac{n+1}{n}l(\sigma)} h_\sigma,\quad
\text{and}\quad 
e_-=\sum_{\sigma\in S_k} (-q^\frac{1-n}{n})^{l(\sigma)} h_\sigma,
\end{equation}
and shows that
\begin{equation}\label{e}
h_\sigma e_+=e_+ h_\sigma= q^{\frac{n-1}{n}l(\sigma)}e_+ \quad \text{and}\quad
h_\sigma e_-=e_- h_\sigma= (-q^{-\frac{n+1}{n}})^{l(\sigma)}e_-,
\end{equation}
for any $\sigma\in S_k.$
Furthermore, for $P_\pm=\sum_{\sigma\in S_k} q^{\pm 2l(\sigma)},$
$e_{\pm}/P_\pm$ are primitive idempotents of $H_k$:
$e_+/P_+$ is the symmetrizer and $e_-/P_-$ is the antisymmetrizer.

%
\section{Proof of Theorem \ref{main}}
\label{s_proof_main}
%

The isotopy invariance of the $n$-bracket follows from 
\cite[Theorem 5.1]{RT} and from the following proposition.

\begin{proposition}\label{tensors}$\phantom{99}$
\begin{itemize} 
\item[\rm(i)] The tensors (5.1.1)-(5.1.3) in \cite{RT} for
the defining $U_h$-representation are given by our formulas 
(\ref{crossings})-(\ref{cupscaps}). (Recall that $q=e^\frac{h_{RT}}{2},$ where
$h_{RT}$ is the Reshetikhin-Turaev $h.$)
\item[\rm(ii)] The map (\ref{T+}), ${\cal T}_+: \C(q)\to V^{\otimes n},$ is 
$U_q$-equivariant.
\item[\rm(iii)]The map (\ref{T-}), ${\cal T}_-: (V^*)^{\otimes n}\to \C(q),$
is $U_q$-equivariant.
\end{itemize}
\end{proposition}

\noindent {\bf Proof of Proposition \ref{tensors}(i)}\qua
By \cite{KS}, the $R$-matrix acts on 
$V\otimes V$ by the matrix 
\begin{equation}\label{R-matrix2}
q^{-\frac{1}{n}}\left(q\sum_i (E_{ii}\otimes E_{ii}) +\sum_{i\ne j} 
(E_{ii}\otimes E_{jj}) + (q-q^{-1})\sum_{i>j} (E_{ij}\otimes E_{ji})\right).
\end{equation}
(This matrix is denoted by $R_{1,1}$ in \cite{KS}, cf.\ the first paragraph
of Section 8.4.2 and (60) in \cite{KS}.)
Here, as before, $E_{ij}$ represents the map
$\delta^j_i: V\to V,$ $\delta^j_i(e_k)=\delta_{j,k} e_i.$ 
By composing the map represented by (\ref{R-matrix2}) with the transposition
$\tau: V\otimes V\to V\otimes V,$ $\tau(v_1,v_2)=(v_2,v_1),$ 
we obtain the map $\hat R$ given by (\ref{R-matrix}).

The ``cap'' maps in \cite[(5.1.1-2)]{RT}
are given by the contraction map $V^*\otimes V\to \C[[h]],$
$(x,y)=x(y)$ and the map $V\otimes V^*\to \C[[h]],$ $(y,x)\to x(v^{-1}uy).$
The ``cup'' maps are their duals.
Let us recall the meaning of $v^{-1}u$ in \cite{RT}:
Let $\rho$ be the half-sum $\frac{1}{2}\sum_{\alpha\in \Delta_+}
\alpha$ of primitive roots of $sl_n$ and let $\rho_i\in \Z$
be the coordinates of $\rho$ in the basis of the Cartan subalgebra of $sl_n$
given by simple roots $\alpha_1,...,\alpha_{n-1}.$ By \cite[(7.1.1)]{RT},
$v^{-1}u=exp(2h\sum_{i=1}^{n-1}\rho_iH_i).$ (Recall that
Reshetikhin-Turaev $h$ is twice the $h$ we use.)
Therefore the following lemma completes the proof of Proposition 
\ref{tensors}(i).

\begin{lemma} $v^{-1}u$ acts on $V$ by sending $e_k$ to
$q^{2k-n-1}e_k$
\end{lemma}

\begin{proof}
Positive roots in the Cartan subalgebra of $sl_n$ are of the form
$\alpha_{ij}=\alpha_i+\alpha_{i+1}+...+\alpha_j,$ for 
$1\leq i\leq j \leq n-1.$ Therefore, $\alpha_k$ appears 
$\frac{1}{2}k(n-k)$ times in $\rho$ and $\rho_k=\frac{1}{2}k(n-k).$

By \cite[Section 8.4.1]{KS}, $H_ie_k$ is $-e_k$ if $i=k,$ $e_k$ if
$i=k-1,$ and $0$ otherwise. Therefore, 
$$exp(hH_i)e_k=\begin{cases} q^{-1}e_k & \text{if i=k}\\
qe_k & \text{if i=k-1}\\
e_k & \text{otherwise}.
\end{cases}$$
Consequently, $(v^{-1}u)e_k=q^{-2\rho_k+2\rho_{k-1}}e_k=
q^{(k-1)(n-k+1)-k(n-k)}e_k=q^{2k-n-1}e_k.$
\end{proof}

\noindent{\bf Proof of Proposition \ref{tensors}(ii) -- $U_q$-equivariance of 
${\cal T}_+$}  

$U_q$ acts on the $n$-th power of the defining representation $V,$ via 
the map $$U_q\stackrel{\Delta^{n-1}}{\longrightarrow} U_q^{\otimes n}
\stackrel{\rho^{\otimes n}}{\longrightarrow} End(V^{\otimes n}),$$
where $\Delta^{n-1}:U_q\to U_q^n$ is the $(n-1)$st power of the 
comultiplication in $U_q.$ The following explicit formulas for 
$\Delta^{n-1}$ follow by induction on $n$ from the 
definition of $\Delta,$ \cite[Prop 6.1.2.5]{KS}: 
\begin{gather}\label{KN}
\Delta^{n-1}(K_i)=K_i\otimes ...\otimes K_i, \text{for $i=1,...,n-1,$}\\
\label{EN}
\Delta^{n-1}(E_i)=\sum_{j=1}^n 1\otimes ...\otimes 1\otimes \stackrel{j}{E_i}
\otimes K_i \otimes ...\otimes K_i,\text{for $i=1,...,n-1,$}
\end{gather}
where the index $j$ over $E_i$ means that it takes the $j$-th position
in the tensor product.
For $\Delta^{n-1}(F_i)$ we have a similar expression:
\begin{equation}\label{FN}
\Delta^{n-1}(F_i)=\sum_{j=1}^n K_i^{-1}\otimes ...\otimes K_i^{-1}
\otimes \stackrel{j}{F_i} \otimes 1 \otimes ...\otimes 1.
\end{equation}
Since $U_q$ acts on $\C(q)$ by counit map, $\epsilon:U_q\to \C(q),$ which
sends $E_i,F_i$ to $0$ and $K_i$ to $1,$ we need to show that
$\Delta^{n-1}(K_i)T_+=T_+,$ and
$\Delta^{n-1}(E_i)T_+=\Delta^{n-1}(F_i)T_+=0$ for $i=1,...,n-1.$
In order to prove the first equality notice that by (\ref{K-action}) and
(\ref{KN}) $\Delta^{n-1}(K_i)$ multiplies the $e_i$ component
in $e_{\sigma(1)}\otimes... \otimes e_{\sigma(n)}$ 
by $q^{-1}$ and it multiplies the $e_{i+1}$ component by $q.$ Since it leaves
all other components unchanged,
$\Delta^{n-1}(K_i)\cdot e_{\sigma(1)}\otimes... \otimes e_{\sigma(n)}=
e_{\sigma(1)}\otimes... \otimes e_{\sigma(n)}$ and, consequently,
$\Delta^{n-1}(K_i)T_+=T_+.$
We complete the proof by showing that $\Delta^{n-1}(E_i)T_+=0.$
The proof of $\Delta^{n-1}(F_i)T_+=0$ is analogous and left to the reader.

For simplicity, denote $e_{\sigma(1)}\otimes... \otimes e_{\sigma(n)}\in
V^{\otimes n}$ by $e_{\sigma}.$ By (\ref{length}),
$$l((i,i+1)\sigma)=\begin{cases}
l(\sigma)+1 & \text{if $\sigma^{-1}(i)<\sigma^{-1}(i+1)$}\\
l(\sigma)-1 & \text{otherwise.}
\end{cases}$$
Therefore,
$$T_+=\sum_{\sigma\in S_n \atop \sigma^{-1}(i)<\sigma^{-1}(i+1)}
(-q)^{l(\sigma)}(e_{\sigma}-qe_{(i,i+1)\sigma}),$$
and our goal is to prove that
\begin{equation}\label{t+last_step}
\Delta^{n-1}(E_i)\cdot (e_{\sigma}-qe_{(i,i+1)\sigma})=0,
\end{equation}
for $\sigma$ such that $\sigma^{-1}(i)<\sigma^{-1}(i+1).$
We have
$$1\otimes ...\otimes 1\otimes \stackrel{j}{E_i} \otimes K_i ...\otimes K_i
\cdot e_\sigma=$$
{\small$$\begin{cases}
0 & \text{if $\sigma(j)\ne i$}\\
e_{\sigma(1)}\otimes...\otimes e_{\sigma(j-1)}\otimes \stackrel{j}{e_{i+1}}
\otimes e_{\sigma(j+1)}\otimes ...
\otimes e_{\sigma(n)} & \text{if $\sigma(j)=i$ and $\sigma^{-1}(i+1)<j$}\\
qe_{\sigma(1)}\otimes...\otimes e_{\sigma(j-1)}\otimes \stackrel{j}{e_{i+1}}
\otimes e_{\sigma(j+1)}\otimes ...
\otimes e_{\sigma(n)} & \text{if $\sigma(j)=i$ and $\sigma^{-1}(i+1)>j.$}\\
\end{cases}$$}%
Therefore, if $\sigma^{-1}(i)=j$ then
$\Delta^{n-1}(E_i)\cdot e_{\sigma}=$
$$e_{\sigma(1)}\otimes... \otimes e_{\sigma(j-1)} \otimes \stackrel{j}{e_{i+1}}
\otimes e_{\sigma(j+1)}\otimes ...
\otimes e_{\sigma(n)} 
\begin{cases}
q & \text{if $\sigma^{-1}(i)<\sigma^{-1}(i+1)$}\\
1 & \text{otherwise.}\\
\end{cases}$$
This implies (\ref{t+last_step}) and, therefore, completes the proof
of $U_q$-equivariance of ${\cal T}_+.$

\medskip
\noindent{\bf Proof of Proposition \ref{tensors}(iii) -- $U_q$-equivariance 
of ${\cal T}_-$}

We need to prove that for any $x\in U_q$ and any
$w=v_{j_1}\otimes ...\otimes v_{j_n}
\in V^{\otimes n},$ $\epsilon(x){\cal T}_-(w)={\cal T}_-(\Delta^{n-1}(x)\cdot 
w).$ This equality reduces to the following three sets of
equations for $i=1,...,n-1:$
\begin{gather}\label{K-eq}
{\cal T}_-(\Delta^{n-1}(K_i)\cdot w)={\cal T}_-(w)\\
\label{E-eq}{\cal T}_-(\Delta^{n-1}(E_i)\cdot w)=0\\
\label{F-eq}{\cal T}_-(\Delta^{n-1}(F_i)\cdot w)=0.
\end{gather}
Both sides of (\ref{K-eq}) vanish if the numbers $(i_1,...,i_n)$ are
not distinct. On the other hand, if these numbers are distinct then
$\Delta^{n-1}(K_i)\cdot w=w$ and (\ref{K-eq}) follows.
We will complete the proof by showing (\ref{E-eq}) -- the proof of 
(\ref{F-eq}) is analogous.

Observe that $\Delta^{n-1}(E_i)\cdot v_{j_1}\otimes ...\otimes v_{j_n}$
is a linear combination of terms $v_{k_1}\otimes ...\otimes v_{k_n}$
such that the $n$-tuple $(k_1,...,k_n)$ is obtained from
$(j_1,...,j_n)$ by changing one of the indices from $i$ to $i+1.$
Since ${\cal T}_-(v_{k_1}\otimes ...\otimes v_{k_n})=0$
if the numbers $k_1,...,k_n$ are not a permutation of 
$1,...,n,$ the left side of (\ref{E-eq}) vanishes unless
$j_1,...,j_{l-1},j_l+1,j_{l+1},...,j_n$ are a permutation $\sigma$ of 
$1,...,n,$ for some $l$ such that $j_l=i.$ In this case $j_k=i$ for
some $k\ne l$ and we can assume that $k<l.$ Under above assumptions, 
\begin{eqnarray*}
\Delta^{n-1}(E_i)\cdot (v_{j_1}\otimes ...\otimes v_{j_n}) & =&
v_{j_1}\otimes ...\otimes \stackrel{k}{v_{i+1}}\otimes ...\otimes 
\stackrel{l}{q^{-1}v_i}\otimes ... \otimes v_{j_n} +\\
& & v_{j_1}\otimes ...\otimes \stackrel{k}{v_i}\otimes ...\otimes 
\stackrel{l}{v_{i+1}} \otimes ...\otimes v_{j_n}.\\
\end{eqnarray*}
Hence, ${\cal T}_-\left(\Delta^{n-1}(E_i)\cdot (v_{j_1}\otimes ...
\otimes v_{j_n})\right)=(-q)^{l((i,i+1)\sigma)}\cdot q^{-1}+
(-q)^{l(\sigma)}.$ Since $l((i,i+1)\sigma)=l(\sigma)+1,$ 
the left hand side of (\ref{E-eq}) vanishes and the proof 
of Proposition \ref{tensors}(iii) is completed.
\medskip

\noindent{\bf Proof of Theorem \ref{main}(ii)}\qua
In the previous section, we proved that $\lb\cdot \rb_n$ is an 
isotopy invariant of
$n$-webs. Now we are going to show that it satisfies properties (i)-(v)
formulated in Theorem \ref{skein_def}.

\noindent{\bf (i)}\qua Since Klimyk's and Schm\"udgen's $\hat R$ is our
$q^\frac{1}{n}\hat R,$ our $\hat R$ 
satisfies 
$$(q^{\frac{1}{n}}\hat R-q)(q^{\frac{1}{n}}\hat R+q^{-1})=0$$
by \cite[Proposition 8.4.24]{KS}
and, hence, 
$$q^\frac{1}{n}\hat R - q^{-\frac{1}{n}}{\hat R}^{-1}=
(q-q^{-1})I.$$ This implies the skein relation (i) of Theorem
\ref{skein_def}.

\noindent{\bf (ii)}\qua Since the two relations (ii) of Theorem \ref{skein_def} 
are inverses of each other, we will show the first of them only.

The ``kink,'' \diag{kink+}{.3in}, defines a map on $V$ 
which is $U_q$-equivariant.
Since $V$ is an irreducible module, this map is a multiple of $Id_V$
and, therefore, for our purpose it is enough to show that the kink maps 
$e_i$ to $q^{n-\frac{1}{n}}e_i$ for some (and hence for arbitrary) $i.$ 
Choose $i=n.$ 
Since the arc \diag{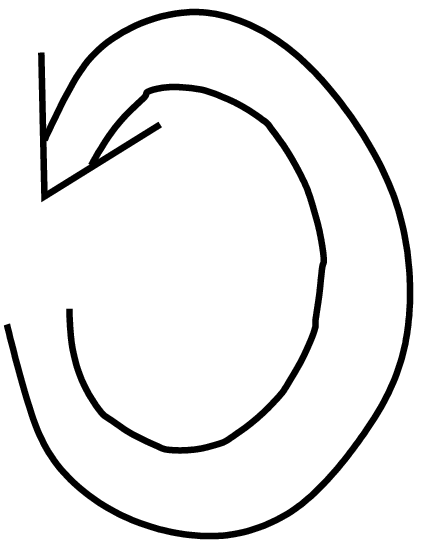}{.3in} maps $V$ to itself by sending $e_i$ to 
$q^{2i-n-1}e_i,$
the kink maps $e_n$ to $Ce_n,$ where
$C=\sum_{k=1}^n {\hat r}_{nk}^{nk} q^{2k-n-1}$ and ${\hat r}_{ij}^{nk}$
are the coefficients of the $\hat R$-matrix,
$${\hat R}(e_n\otimes e_k)=
\sum_{i,j} {\hat r}_{ij}^{nk}e_i\otimes e_j.$$
Since ${\hat r}^{nk}_{nk}=\begin{cases}
q^{-\frac{1}{n}+1} & \text{for $k=n$}\\
0 & \text{otherwise,}\\
\end{cases}$
$C=q^{n-\frac{1}{n}}.$

\noindent {\bf (iii)}\qua This property will be proved in the next
section.

{\bf (iv)}\qua 
The  bracket of the trivial knot diagram is given by the
contraction of
the cup and the cap tensors, where the cup and the cap are chosen with 
coinciding orientations. Therefore 
$$\lb \diag{annulus}{.2in}\rb_n= \sum_{i=1}^n q^{2i-n-1}=[n].$$
By the construction of the bracket, $\lb \Gamma_1 \cup \Gamma_2\rb_n=
\lb \Gamma_1\rb_n \cdot \lb \Gamma_2\rb_n,$
for disjoint (and hence unlinked) web diagrams $\Gamma_1,\Gamma_2.$

{\bf (v)}\qua This is obvious.
\endproof

%
\section{Proof of Proposition \ref{extra_rel} and of 
Theorem \ref{skein_def}(iii)}\label{main_proofs}
\label{p_skein_def}
%

We will often use the following equality
\begin{equation}\label{gy}
\sum_{\sigma\in S_n} (-q)^{2l(\sigma)}=
q^\frac{n(n-1)}{2}\cdot [n]!
\end{equation}
following from \cite[(3.1)]{Gy}.

\begin{lemma}\label{basicweb}
$\lb \diag{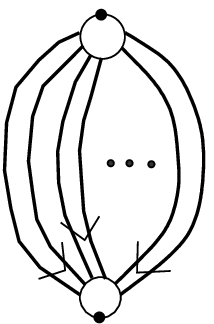}{.6in}\rb_n=q^\frac{n(n-1)}{2}\cdot [n]!$
\end{lemma}

\begin{proof}
The above bracket is given by the contraction of ${\cal T}_-$ with
${\cal T}_+,$ 
$${\cal T}_-(T_+)=\sum_{\sigma\in S_n} (-q)^{2l(\sigma)}=
q^\frac{n(n-1)}{2}\cdot [n]!,$$
by (\ref{gy}).
\end{proof}

Consider the skein 
\begin{equation}\label{antisym}
\Lambda_k=\sum_{\sigma\in S_k} 
(-q^\frac{1-n}{n})^{l(\sigma)}\diag{sigma}{.4in},
\end{equation}
where, as before, \diag{sigma}{.4in} is the unique positive braid with
$l(\sigma)$ crossings representing $\sigma.$
Given an $(k,k)$-tangle $T,$ \diag{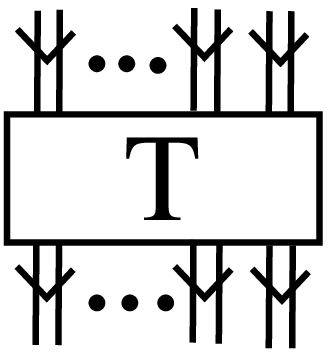}{.4in}, denote
by $\pi_k(T)$ the \mbox{$(k-1,k-1)$}-tangle obtained from $T$ by closing 
up its last string, \diag{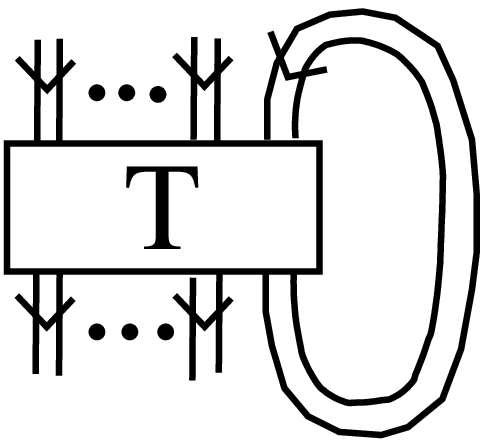}{.5in}\, . The definition of 
$\pi_k(T)$ obviously extends to all skeins $T$ being linear combinations 
of $(k,k)$-tangles.

\begin{lemma}\label{closure} $\pi_{k+1}(\Lambda_{k+1})=\Lambda_k q^{-k}[n-k].$
\end{lemma}

\begin{proof}
Each permutation $\sigma\in S_{k+1}$ can be written in the form
\begin{equation}\label{permutationform}
(i_1,i_1-1,...,i_1-j_1)(i_2,i_2-1,...,i_2-j_2)...(i_l,i_l-1,...,i_l-j_l),
\end{equation}
where $i_1<i_2<...<i_l.$ Furthermore, such presentation is unique.
(These statements can be proved by induction on $k$.) 
By splitting the set of all permutations $\sigma\in S_{k+1}$ into those 
with $i_l\leq k$ and those with $i_l=k+1,$ we get
$$\pi_{k+1}(\Lambda_{k+1})=\Lambda_k [n]+
\sum_{i=1}^k (-q^{\frac{1-n}{n}})^{k+1-i}\diag{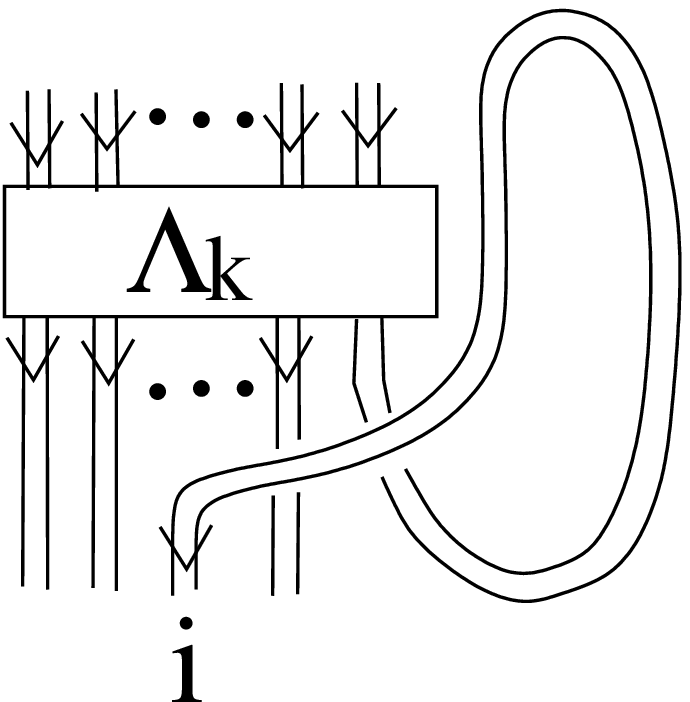}{.8in}.$$
The $i$-th summand in the sum on the right side takes into account all 
permutations $\sigma\in S_{k+1}$ with $i_l=k+1$ and $j_l=i.$
Note that the action of $\Lambda_k$ on $V^{\otimes k}$ coincides with the
one of $e_-,$ defined in (\ref{e-def}), and hence composing $\Lambda_k$
with a single positive crossing on two adjacent strings yields
$-q^{-\frac{n+1}{n}}\Lambda_k.$ 
Therefore, after applying relation (ii) of Theorem \ref{skein_def}
to remove the kink in the skein above and after replacing the $k-i$ 
crossings by the factor $(-q^{-\frac{n+1}{n}})^{k-i}$ we get 
$$\pi_{k+1}(\Lambda_{k+1})=\Lambda_k [n]+
\sum_{i=1}^k (-q^\frac{1-n}{n})^{k+1-i} (-q^{-\frac{n+1}{n}})^{k-i}
q^{n-n^{-1}} \Lambda_k=$$
$$\Lambda_k [n]+ \Lambda_k q^{n-n^{-1}} (-q^\frac{1-n}{n})\sum_{i=1}^k 
q^{-2(k-i)}=\Lambda_k [n]+ \Lambda_k (-q^{n-1})\frac{1-q^{-2k}}{1-q^{-2}}=$$
$$\Lambda_k \left([n]-\frac{q^n-q^{n-2k}}{q-q^{-1}}\right)=
\Lambda_k \left(\frac{q^{n-2k}-q^{-n}}{q-q^{-1}}\right)=
\Lambda_k q^{-k}[n-k].$$
\end{proof}

Recall that $T_+\in V^{\otimes n}$ was defined in (\ref{T+})
and $h_{i,i+1}$ in Section \ref{s_hecke}.

\begin{corollary}\label{antisym_closure}
$$\sum_{\sigma\in S_n} 
(-q^\frac{1-n}{n})^{l(\sigma)}\lb \diag{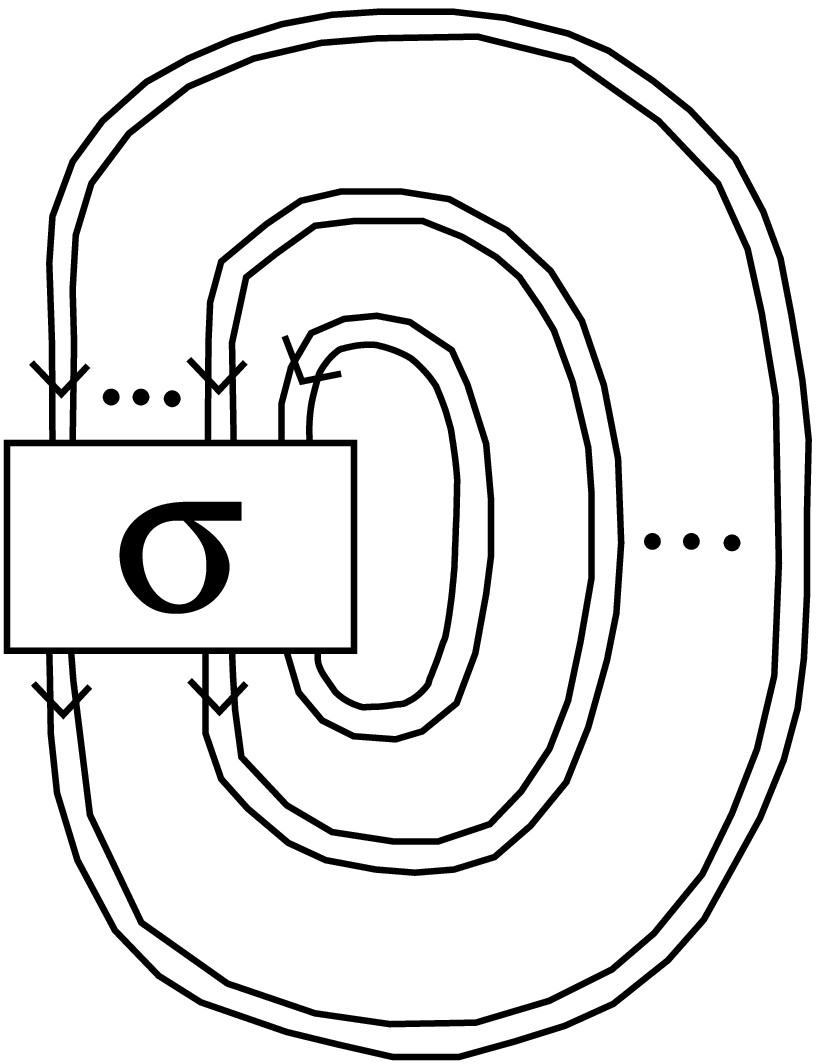}{.7in}\rb_n=
q^{-\frac{n(n-1)}{2}}[n]!$$
\end{corollary}

\begin{lemma}\label{hT_+}
$h_{(i,i+1)}T_+=(-q^{-\frac{n+1}{n}})T_+.$
\end{lemma}

\begin{proof}
Denote $e_{\sigma(1)}\otimes ... e_{\sigma(n)}\in V^{\otimes n}$
by $e_\sigma$ as before. By Proposition \ref{prop_hecke}(1) and
(\ref{R-matrix}), 
$$h_{(i,i+1)}e_\sigma=q^{-\frac{1}{n}}e_{(\sigma(i),\sigma(i+1))\sigma}+
\begin{cases} q^{-\frac{1}{n}}(q-q^{-1})e_{\sigma} & \text{if 
$\sigma(i)<\sigma(i+1)$}\\
0 & \text{otherwise}\\
\end{cases}$$
Let $A_1$ be the set of these
permutations $\sigma\in S_n$ for which $\sigma(i)<\sigma(i+1),$
and let $A_2=S_n\setminus A_1.$
Furthermore, let $$T_i=\sum_{\sigma\in A_i} (-q)^{l(\sigma)}e_\sigma,$$
for $i=1,2.$
Then 
\begin{equation}\label{hT_2}
\begin{array}{c}
h_{(i,i+1)}T_2=q^{-\frac{1}{n}} \sum_{\sigma\in A_2} 
(-q)^{l(\sigma)}e_{(\sigma(i),\sigma(i+1))\sigma}=\\
q^{-\frac{1}{n}} \sum_{\sigma\in A_2} 
(-q)^{l((\sigma(i),\sigma(i+1))\sigma)+1}e_{(\sigma(i),\sigma(i+1))\sigma}=
q^{-\frac{1}{n}}(-q)T_1.
\end{array}
\end{equation}
Similarly,
\begin{equation}\label{hT_1}
h_{(i,i+1)}T_1=q^{-\frac{1}{n}} \sum_{\sigma\in A_1} 
(-q)^{l(\sigma)}e_{(\sigma(i),\sigma(i+1))\sigma}+q^{-\frac{1}{n}}(q-q^{-1})
\sum_{\sigma\in A_1} (-q)^{l(\sigma)}e_{\sigma}.
\end{equation}
Note that $\sigma\in A_1 \Leftrightarrow (\sigma(i),\sigma(i+1))\sigma\in A_2$
and $l((\sigma(i),\sigma(i+1))\sigma)=l(\sigma)+1,$ for $\sigma\in
A_1.$ 
Therefore,
after substituting $\tau=(\sigma(i),\sigma(i+1))\sigma$ in the first
sum of (\ref{hT_1}) we get
\begin{eqnarray*}
h_{(i,i+1)}T_1 & =& q^{-\frac{1}{n}} \sum_{\tau\in A_2} 
(-q)^{l(\tau)-1}e_{\tau}+q^{-\frac{1}{n}}(q-q^{-1}) T_1\\
 & =& q^{-\frac{1}{n}}(-q^{-1})T_2+ q^{-\frac{1}{n}}(q-q^{-1}) T_1.\\
\end{eqnarray*}
Hence, by (\ref{hT_2}), $h_{(i,i+1)}(T_1+T_2)=(-q^{-\frac{n+1}{n}})(T_1+T_2).$
\end{proof}

\noindent{\bf Proof of Theorem \ref{skein_def}(iii)}\qua
We need to prove that the skeins $S= \diag{sinksource}{.5in}$ and
$q^{n(n-1)}\Lambda_n$ 
coincide as operators on $V^{\otimes n}.$
Since $S$ is $U_q$-equivariant, by Proposition
\ref{prop_hecke}(iii), it is is equal to the map $w\to x\cdot w: 
V^{\otimes n}\to V^{\otimes n},$ for certain $x\in H_n.$ 
Since the image of $S$ is $1$-dimensional, $x$ is either a multiple
of $e_+$ or $e_-.$ (This statement follows from the fact that
$H_n$ is isomorphic to the group ring of $S_n$ over $\C(t).$)
Lemma \ref{hT_+} indicates that $S$ is a multiple of $e_-,$ ie.
$S=ce_-$ for certain $c\in \C(q^\frac{1}{n}).$
On the other hand, $\Lambda_n$ coincides with $e_-$ as an operator on 
$V^{\otimes n}.$
Therefore, we need to prove that 
$c=q^{n(n-1)},$ and for that it is enough to consider the closures of 
$S$ and $\Lambda_n.$
Now, the statement follows from Lemma \ref{basicweb} and 
Corollary \ref{antisym_closure}.
\endproof

\noindent {\bf Proof of Proposition \ref{extra_rel}}\qua
By Theorem \ref{skein_def}(iii) and Lemma \ref{closure}, 
$$\diag{crossres}{.5in}=q^{n(n-1)}\cdot \pi_n(... \pi_3(\Lambda_n))=
q^{n(n-1)} q^{-\frac{n(n-1)}{2}+1}[n-2]!\Lambda_2.$$
By substituting $\diag{l0}{.3in}-q^\frac{1-n}{n}\diag{l+}{.3in}$ for 
$\Lambda_2$ we get the statement of Proposition \ref{extra_rel}.
\endproof

%
\section{Proof of Proposition \ref{marked}}
\label{s_marked_proof}
%

We prove the statement for sources only -- the proof for sinks is
analogous. We begin with two preliminary results.

Recall that for any two sets of integers, $S_1,S_2,$  
$$\pi(S_1,S_2)=\#\{(i_1,i_2)\in S_1\times S_2: i_1>i_2\}.$$
The proof of the following lemma is left to the reader:

\begin{lemma}\label{pilemma}
If $S=\{s_1,...,s_k\}\subset N=\{1,...,n\}$ and $S'=N\setminus S,$
then $\pi(S,S')=\sum_{i=1}^k s_i - \frac{k(k+1)}{2}.$
\end{lemma}

Let $$\tau_{n,k}=\left(\begin{array}{cccccc}
1&... & k,& k+1& ...& n\\
n-k+1&...&n,&1&...& n-k\\
\end{array}\right).$$

\begin{lemma}\label{lsigma} For any $\sigma\in S_n,$
\begin{equation}
l(\sigma)=l(\sigma\tau_{n,k})+ 2\sum_{i=1}^k \sigma(i)-k(n+1).
\end{equation}
\end{lemma}

\begin{proof} Fix $k<n.$
By (\ref{length}), $l(\sigma)=A+B+C,$
where 
$$A=\#\{i<j\leq n-k: \sigma(i)>\sigma(j)\},\ 
B=\#\{n-k< i<j\leq n: \sigma(i)>\sigma(j)\},$$ 
$$C=\#\{i\leq n-k<j: \sigma(i)>\sigma(j)\}.$$
Similarly, 
$$l(\sigma\tau_{n,k})=\#\{i<j\leq k: \sigma\tau_{n,k}(i)>
\sigma\tau_{n,k}(j)\}+$$
$$\#\{k< i<j\leq n: \sigma\tau_{n,k}(i)>\sigma\tau_{n,k}(j)\} + 
\#\{i\leq k< j: \sigma\tau_{n,k}(i)>\sigma\tau_{n,k}(j)\}.$$
Note that the first, second, and the third summands above are equal $B,$ 
$A,$ and $k(n-k)-C$ respectively.
By Lemma \ref{pilemma}, $C=\sum_{i=1}^k \sigma(i)-\frac{k(k+1)}{2},$ 
and hence the statement follows.
\end{proof}

Let ${\cal T}_{-}':V^{\otimes n}\to \Z(q)$ be the tensor
associated with \diag{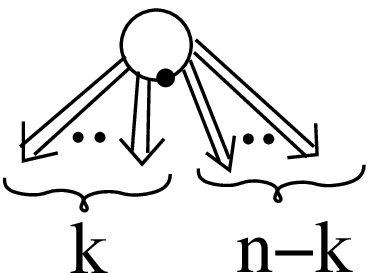}{.5in}

for some $k.$
We need to show that ${\cal T}_-'={\cal T}_-$ for $n$ odd and
${\cal T}_-'={\cal T}_-$ mod $2$ for $n$ even, where
${\cal T}_-:V^{\otimes n}\to \Z(q)$ is the tensor defined 
by (\ref{sourcesink}) in Section \ref{sqdefinition}.
Here is another presentation of the above graph:
\diag{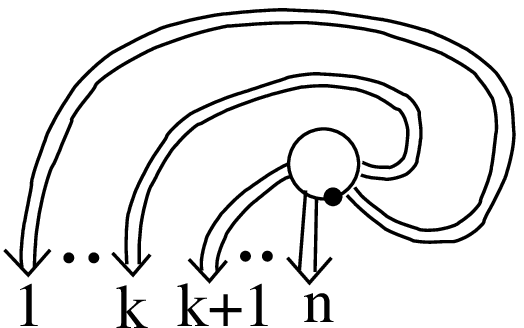}{.6in}
 
The tensor ${\cal T}_-'$ is given by the contraction of cups and caps
placed on strings $1,...,k$ with the tensor ${\cal T}_-.$
Hence ${\cal T}_-':V^{\otimes n}\to \Z(q)$ equals 
${\cal T}_-\Psi,$ where $$\Psi(e_\sigma)=e_{\sigma\tau_{n,k}}\cdot
q^{2\sum_{i=1}^k \sigma(i)-k(n+1)}.$$
Since ${\cal T}_-(e_\sigma)=(-q)^{l(\sigma)}$ and
$${\cal T}_-(\Psi(e_\sigma))=(-q)^{l(\sigma\tau_{n,k})}
q^{2\sum_{i=1}^k \sigma(i)-k(n+1)}.$$
Now the statement follows from Lemma
\ref{lsigma}.
\endproof

%
\section{Proofs of Lemma \ref{core_indep} and Theorem \ref{statesum}:}
\label{p_statesum}
%

For any core $\Gamma'$ of an $n$-web $\Gamma$ let 
\begin{equation}\label{the_sum}
rot_{n,\Gamma'}(S)=\sum_e ind(e)(2S(b)-n-1)\in \Z.
\end{equation}

\begin{lemma}
For any core $\Gamma'$ of $\Gamma,$ 
$rot_{n,\Gamma'}(S)\in \Z.$
\end{lemma}

\begin{proof}
Let $\Gamma'$ be a core of $\Gamma$ and let $v$ be the marked
point of a disc of $\Gamma.$ Suppose that $\Gamma$ is isotoped to
$\wt\Gamma$ and $\Gamma'$ is isotoped to a core $\wt\Gamma'$ of
$\wt\Gamma$ such that the tangents at the endpoints of edges of
$\wt\Gamma'$ are unchanged, except for those at $v.$ Then the indices, 
$ind(e),$ remain unchanged, except for those edges $e$ which are 
adjacent to $v.$ For these edges
$ind_{\wt\Gamma'}(e)=ind_{\Gamma'}(e)+ \beta,$ for
some $\beta$ (which is the same for all edges $e$ adjacent to $v$).
Since $$\sum_{\text{$e$ adjacent to $v$}} \beta (2S(b)-n-1)=0,$$
$$rot_{n,\wt \Gamma'}(S)=rot_{n,\Gamma'}(S).$$
Fix a vector $\vec w.$ By performing appropriate isotopies of $\Gamma$ 
and of $\Gamma'$ we may assume that $\Gamma'$ is such that the tangent 
vector to any
endpoint $v$ of every edge of $\Gamma'$ is either $\vec w$ or
$-\vec w$ depending if $v$ is the marked point of a source or a sink.
In this situation, $ind(e)\in \Z$ for all edges $e$ of $\Gamma'$
and, consequently, 
$$rot_{n,\Gamma'}(S)\in \Z.$$
\end{proof}

Since any two cores of $\Gamma$ are isotopic to each other and
$rot_{n,\Gamma'}(S)$ varies continuously under isotopy of
$\Gamma',$ $rot_{n,\Gamma'}(S)$ is independent of the
choice of $\Gamma'.$ This completes the proof
of Lemma \ref{core_indep}.\vspace*{.1in}

\noindent{\bf Proof of Theorem \ref{statesum}}\qua
By Theorem \ref{main}(i), $\lb \Gamma\rb_n$ is given by a contraction
of tensors corresponding to the vertices and ``caps'' and ``cups'' of
$\Gamma.$ (By assumption of Theorem \ref{statesum}, $\Gamma$ has no
crossings). Note that the summands in that sum are in 1-1
correspondence with the states of $\Gamma$ and that each of the summands is
a power of $\pm q.$ By deforming $\Gamma$ by an isotopy if
necessary, we can choose a core $\Gamma'$ of $\Gamma$ so that
the tangent vectors to the edges of $\Gamma'$ at their endpoints point all
in the same direction. It is easy to see that for such $\Gamma',$ 
$q^{rot_n(S)}$ is the power of $q$ given by the cups and
caps of $\Gamma.$  Furthermore, any state $S$ and any vertex $v,$ 
the tensor associated with $v$ contributes $(-q)^{l(P(S,v))}$ 
to the state sum.
\endproof

%
\section{Proof of Proposition \ref{MOY-state_sum}}
\label{sMOY-state_sum_proof}
%

For any $k>0$ we identify the basis vectors $e_{i_1}\otimes ...\otimes
e_{i_k}\in V^{\otimes k}$ with sequences $(i_1,...,i_k)\in
\{1,...,n\}^k.$ For any $a=(a_1,...,a_k)$ we denote the
set $\{a_1,...,a_k\}$ by $\bar a.$ Furthermore, we denote 
$\{1,...,n\}$ by $\bar n.$

An {\em enhanced state} $S$ of $MOY_n$-graph $\Gamma$ is a function 
which assigns to each edge $e$ a sequence $S(e)=(i_1,...,i_{|e|})$
of $|e|$ distinct elements of the set $\bar n,$
such that $\bar S(e_1^v)\cup \bar S(e_2^v)= \bar S(e_0^v)$ for any vertex $v.$
Any enhanced state $S$ defines a state $\bar S$ of $\Gamma$ 
labeling every edge $e$ of $\Gamma$ by the set $\bar S(e).$ 

For a sequence $a$ of numbers $(a_1,...,a_k),$ which does not
contain any repetitions, we denote by $l(a)$ the
length of the permutation which puts the numbers of the sequence in
the increasing order. 

Denote the tensors associated with the graphs
\begin{equation}\label{T1T2}
\diag{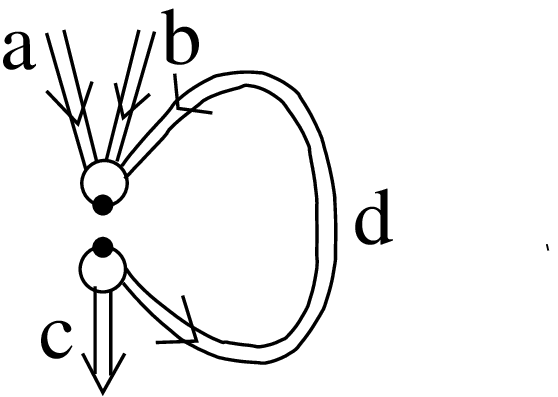}{.7in}\quad \quad \diag{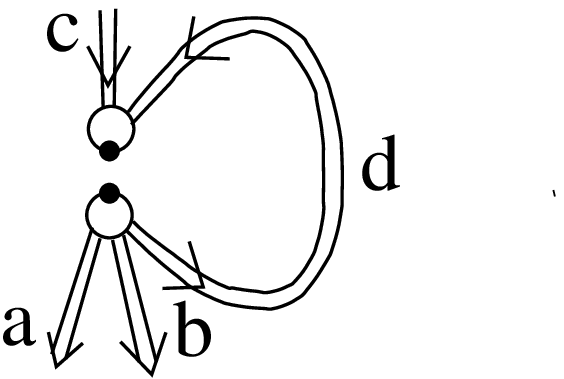}{.65in}
\end{equation}
by $T_1: V^{\otimes k+l}\to V^{\otimes k}\otimes 
V^{\otimes l},$ $T_2: V^{\otimes k}\otimes 
V^{\otimes l}\to  V^{\otimes k+l}.$
There are obtained by a partial contraction of
four tensors: a cap, a cup, a sink, and a source,
and are given by the following formulas
$$T_1(c)=\sum_{a,b,d} t(a,b,c,d) a\otimes b,\quad
T_2(a,b)=\sum_{c,d} t(a,b,c,d) c,$$
where $a\in {\bar n}^k,$ $b\in {\bar n}^l,$ $c\in {\bar n}^{k+l},$
$d\in {\bar n}^m,$ $m=n-k-l$ and $t(a,b,c,d)$ is defined as follows:
$t(a,b,c,d)=0$ unless

\begin{itemize}
\item[(i)] $\bar a\cup \bar b \cup \bar d =\bar n$
\item[(ii)] $\bar c=\bar a\cup \bar b.$
\end{itemize}

(The first condition implies that
$\bar a, \bar b,$ and $\bar d$ are disjoint and the sequences
$a,b,d$ have no repeating elements.)
If these conditions are satisfied then
\begin{eqnarray*}
\lefteqn{t(a,b,c,d)=}\\
 & & (-q)^{l(a)+l(b)+l(d)+\pi(\bar a,\bar b)+
\pi(\bar a\cup \bar b,\bar d)}\cdot (-q)^{l(c)+l(d)+\pi(\bar c,\bar
  d)}\cdot q^{2\sum_{i=1}^m d_i- (n+1)m}.
\end{eqnarray*}
The first two factors above come from the tensors associated with 
the vertices of graphs (\ref{T1T2}). The third factor corresponds
to caps and cups of graphs (\ref{T1T2}).
And since $\pi(\bar a\cup \bar b,\bar d)=\pi(\bar c,\bar d),$
we get
\begin{equation}\label{t1}
t(a,b,c,d)=(-q)^{l(a)+l(b)+l(c)+\pi(\bar a,\bar b)}
\cdot q^{2l(d)+2\pi(\bar c,\bar d)+2\sum_{i=1}^m d_i- (n+1)m}
\end{equation}
If conditions (i) and (ii) above are satisfied then by Lemma \ref{pilemma},
$$2(\bar c, \bar d)+2\sum_{i=1}^m d_i - (n+1)m=
2\sum_{i=1}^{k+l} c_i - (k+l)(k+l+1)+
2\sum_{i=1}^m d_i - (n+1)m=$$ 
$$=(n+1)n- (k+l)(k+l+1) -(n+1)m= (n-k-l)(k+l).$$
Hence, under conditions (i) and (ii) above,
$$t(a,b,c,d)=(-q)^{l(a)+l(b)+l(c)+\pi(\bar a,\bar b)}\cdot 
q^{2l(d)+(n-k-l)(k+l)}.$$
By (\ref{gy}), $\sum_{\sigma\in S_m} q^{2l(\sigma)}=
q^\frac{m(m-1)}{2}[m]!$
Therefore, if we denote $\sum_d t(a,b,c,d)$ by $t(a,b,c),$ then
$$t(a,b,c)= (-q)^{l(a)+l(b)+l(c)+\pi(\bar a,\bar b)}\cdot 
q^{\frac{m(m-1)}{2}+(n-k-l)(k+l)}[n-k-l]!=$$
$$=(-q)^{l(a)+l(b)+l(c)+\pi(\bar a,\bar b)}\cdot 
q^\frac{(n-k-l)(n+k+l-1)}{2}[n-k-l]!$$
The complete contraction of tensors associated with
vertices, cups, and caps of $W(\Gamma)$ produces
$$[\Gamma]_n=\lb W(\Gamma)\rb_n=\sum_S \Psi(S),$$
where 
$$\Psi(S)=q^{rot(\bar S)}\cdot \prod_v t(S(e^v_1),S(e^v_2),S(e^v_0))$$
(The first of the above factors is the contraction of the tensors 
associated with ``caps'' and ``cups'' in $\Gamma.$)
Hence
$$\Psi(S)=q^{rot(\bar S)}\cdot (-q)^{\sum_e 2l(S(e)) +\sum_v \pi(\bar S(e_1^v),
  \bar S(e_2^v))} \prod_v q^\frac{(n-|e^v_0|)(n+|e_0^v|-1)}{2}[n-|e_0^v|]!
$$
where the sum $\sum_e$ is taken over all edges of
$\Gamma$ which are not annuli.
For any state $s$ of $\Gamma,$ we denote sum $\sum \Psi(S)$ over 
all enhanced states $S$ such that $\bar S=s,$ by $\Phi(s).$
Hence, by (\ref{gy}),
$$\Phi(s)={\cal N}(\Gamma)\cdot q^{rot(s)} \cdot (-q)^{\sum_v \pi(s(e_1^v),
s(e_2^v))},$$
where $${\cal N}(\Gamma)=
\prod_e q^{\frac{|e|(|e|-1)}{2}}[|e|]!\cdot
\prod_v q^\frac{(n-|e^v_0|)(n+|e_0^v|-1)}{2}[n-|e_0^v|]!$$
Since $$\frac{1}{2}(n-|e^v_0|)(n+|e_0^v|-1)=
\frac{n(n-1)+|e^v_0|-|e^v_0|^2}{2}=$$
$$\frac{1}{2}n(n-1)+\frac{1}{4}\left(|e^v_0|+|e^v_1|+|e^v_2|-
|e^v_0|^2-|e^v_1|^2-|e^v_2|^2\right)-\frac{1}{2}|e^v_1|\cdot|e^v_2|,$$ 
and each edge appears as $e^v_i$ for two different vertices $v,$
$${\cal N}(\Gamma)=
\prod_e q^{\frac{|e|^2-|e|}{2}}[|e|]!\cdot
\prod_e q^\frac{|e|-|e|^2}{2}
\prod_v q^\frac{n(n-1)-|e^v_1|\cdot|e^v_2|}{2}
[n-|e_0^v|]!$$
Hence,
$${\cal N}(\Gamma)=
\prod_e [|e|]!\cdot
\prod_v q^\frac{n(n-1)-|e^v_1|\cdot|e^v_2|}{2}
[n-|e_0^v|]!$$
\endproof

%
\section{Proof of Proposition \ref{MOY-bracket}}
\label{sMOY_bracket_proof}
%

For any $$s: \{\text{edges of $\Gamma$}\}\to 
\text{subsets of $\{1,...,n\}$}$$
let $$\hat s: \{\text{edges of $\Gamma$}\}\to \text{subsets of 
$\{-\frac{n-1}{2},-\frac{n-1}{2}+1,..., \frac{n-1}{2}\}$},$$
be such that $\hat s(e)=\{i_1-\frac{n+1}{2},..., i_k-\frac{n+1}{2}\}$
if $s(e)=\{i_1,...,i_k\}.$ We observed in Section \ref{s_MOY} already,
that $s$ is a state of $\Gamma$ if and
only if $\hat s$ is a MOY-state of $\Gamma.$

Let $\Gamma$ be an $MOY_n$-graph diagram with no crossings.
We can assume that $\Gamma$ is composed of caps, cups, and 
vertices of the following form:
$$\diag{moyine}{.5in}\quad \diag{moyoute}{.5in}$$
Since these pictures are obtained by rotating the pictures of
\cite[Fig. 1.3]{MOY} by $180^o$ and exchanging $e_1$ with $e_2,$
the Murakami-Ohtsuki-Yamada weight associated to the vertices
above is 
$$q^{|e_1^v|\cdot |e_2^v|/4-\pi(s(e_2^v),s(e_1^v))/2}.$$
But since $$\pi(s(e_2),s(e_1))+\pi(s(e_1),s(e_2))=|e_1|\cdot|e_2|,$$
this weight equals to
$$q^{-|e_1|\cdot |e_2|/4+\pi(s(e_1),s(e_2))/2}.$$
Since the Murakami-Ohtsuki-Yamada rotation index of $\hat s$ is
$rot(s),$
$$\{\Gamma\}_n= \sum_{\text{states} s} q^{rot(s)}
\prod_{v} q^{-|e_1^v|\cdot |e_2^v|/4+\pi(s(e_1^v),s(e_2^v))/2},$$
and the statement follows.

\Addressesr


\begin{thebibliography}{99}

\let\olditem\bibitem
\def\bibitem#1{\olditem[#1]{#1}}



\bibitem{APS}
\textbf{M\,M Asaeda}, \textbf{J\,H Przytycki}, \textbf{A\,S
  Sikora}, \emph{Categorification of the {K}auffman bracket skein module of
  {$I$}-bundles over surfaces}, \agtref4{2004}{52}{1177}{1210}
  \MR{MR2113902}

\bibitem{BN}
\textbf{D Bar-Natan}, \emph{On {K}hovanov's categorification of the {J}ones
  polynomial}, \agtref2{2002}{16}{337}{370}
  \MR{1917056}

\bibitem{Bl}
\textbf{C Blanchet}, \emph{Hecke algebras, modular categories and
  {$3$}-manifolds quantum invariants}, Topology 39 (2000) 193--223
  \MR{1710999}

\bibitem{Bu}
\textbf{D Bullock}, \emph{Rings of {${\rm SL}\sb 2({\bf C})$}-characters and
  the {K}auffman bracket skein module}, Comment. Math. Helv. 72 (1997)
  521--542\relax \MR{1600138}

\bibitem{BFK}
\textbf{D Bullock}, \textbf{C Frohman}, \textbf{Joanna
  Kania-Bartoszy{\'n}ska}, \emph{Understanding the {K}auffman bracket skein
  module}, J. Knot Theory Ramifications 8 (1999) 265--277\relax \MR{1691437}

\bibitem{CP}
\textbf{V Chari}, \textbf{A Pressley}, \emph{A guide to quantum
  groups}, Cambridge University Press, Cambridge (1994) \MR{1300632}

\bibitem{FZ} \textbf{C Frohman}, \textbf{J Zhong}, \emph{The
  Yang-Mills measure in the $SU_3$ skein module}, preprint (2004)

\bibitem{FGL}
\textbf{C Frohman}, \textbf{R Gelca}, \textbf{WLofaro},
  \emph{The {A}-polynomial from the noncommutative viewpoint}, Trans. Amer.
  Math. Soc. 354 (2002) 735--747  \MR{1862565}

\bibitem{Ga} \textbf{S Garoufalidis}, Difference and differential
  equations for the colored Jones function, \arxiv{math.GT/0306229}

\bibitem{GL} \textbf{S Garoufalidis}, \textbf{T\,T\,Q Le}, The
  colored Jones function is q-holonomic,
  \gtref{9}{2005}{29}{1253}{1293}


\bibitem{Ge}
\textbf{R Gelca}, \emph{On the relation between the {$A$}-polynomial
  and the {J}ones polynomial}, Proc. Amer. Math. Soc. 130 (2002) 1235--1241
  (electronic)\relax \MR{1873802}



\bibitem{Go} \textbf{B Gornik}, \emph{Note on Khovanov link cohomology},
  \arxiv{math.QA/0402266}

\bibitem{Gy}
\textbf{A Gyoja}, \emph{A {$q$}-analogue of {Y}oung symmetrizer}, Osaka
  J. Math. 23 (1986) 841--852\relax \MR{873212}

\bibitem{HK} \textbf{R\,S Huerfano}, \textbf{M Khovanov}, 
\emph{Categorification of some level two representations of sl(n)}, 
\arxiv{math.QA/0204333}

\bibitem{Ja} \textbf{M Jacobsson}, \emph{An invariant of link cobordisms from 
Khovanov's homology theory}, \agtref4{2004}{53}{1211}{1251} \MR{2113903}

\bibitem{Ka}
\textbf{L\,H Kauffman}, \emph{State models and the {J}ones polynomial},
  Topology 26 (1987) 395--407 \MR{899057}

\bibitem{KV}
\textbf{L\,H Kauffman}, \textbf{P Vogel}, \emph{Link polynomials and a
  graphical calculus}, J. Knot Theory Ramifications 1 (1992) 59--104
  \MR{1155094}

\bibitem{K1}
\textbf{M Khovanov}, \emph{A categorification of the {J}ones polynomial},
  Duke Math. J. 101 (2000) 359--426\relax \MR{1740682}

\bibitem{K2}
\textbf{M Khovanov}, \emph{sl(3) link homology}, 
\agtref4{2004}{45}{1045}{1081}\relax \MR{2100691}

\bibitem{KR} \textbf{M Khovanov}, \textbf{L Rozansky}, \emph{Matrix
  factorizations and link homology}, \arxiv{math.QA/0401268}

\bibitem{KS}
\textbf{A Klimyk}, \textbf{K Schm{\"u}dgen}, \emph{Quantum groups
  and their representations}, Texts and Monographs in Physics, Springer-Verlag,
  Berlin (1997) \MR{1492989}

\bibitem{Ku}
\textbf{G Kuperberg}, \emph{Spiders for rank {$2$} {L}ie algebras}, Comm.
  Math. Phys. 180 (1996) 109--151\relax \MR{1403861}

\bibitem{Le} \textbf{E\,S Lee}, \emph{On Khovanov invariant for alternating
  links}, \arxiv{math.GT/0210213}

\bibitem{LM}
\textbf{A Lubotzky}, \textbf{A\,R Magid}, \emph{Varieties of
  representations of finitely generated groups}, Mem. Amer. Math. Soc. 58
  (1985) xi+117\relax \MR{818915}

\bibitem{Mu} \textbf{H Murakami}, \emph{A quantum introduction to Knot Theory},
  preprint (2003)

\bibitem{MOY}
\textbf{H Murakami}, \textbf{T Ohtsuki}, \textbf{S Yamada},
  \emph{Homfly polynomial via an invariant of colored plane graphs}, Enseign.
  Math. (2) 44 (1998) 325--360\relax \MR{1659228}

\bibitem{OY}
\textbf{T Ohtsuki}, \textbf{S Yamada}, \emph{Quantum {${\rm SU}(3)$}
  invariant of {$3$}-manifolds via linear skein theory}, J. Knot Theory
  Ramifications 6 (1997) 373--404 \MR{1457194}

\bibitem{Pr}
\textbf{J\,H Przytycki}, \emph{Fundamentals of {K}auffman bracket skein
  modules}, Kobe J. Math. 16 (1999) 45--66 \MR{1723531}

\bibitem{PS}
\textbf{J\,H Przytycki}, \textbf{A\,S Sikora}, \emph{On skein algebras
  and {${\rm Sl}\sb 2({\bf C})$}-character varieties}, Topology 39 (2000)
  115--148 \MR{1710996}

\bibitem{Ra} \textbf{J\,A Rasmussen}, \emph{Khovanov homology and the slice
  genus}, \arxiv{math.GT/0402131}

\bibitem{RT}
\textbf{N\,Yu Reshetikhin}, \textbf{V\,G Turaev}, \emph{Ribbon graphs and their
  invariants derived from quantum groups}, Comm. Math. Phys. 127 (1990)
  1--26\relax \MR{1036112}

\bibitem{S1}
\textbf{A\,S Sikora}, \emph{{${\rm SL}\sb n$}-character varieties as spaces
  of graphs}, Trans. Amer. Math. Soc. 353 (2001) 2773--2804 (electronic)\relax
  \MR{1828473}

\bibitem{S2}
\textbf{A\,S Sikora}, \emph{Skein modules and {TQFT}}, from: ``Knots in
  Hellas '98 (Delphi)'', Ser. Knots Everything 24, World Sci. Publishing, River
  Edge, NJ (2000)  436--439\relax \MR{1865721}


\bibitem{Tu}
\textbf{V\,G Turaev}, \emph{The {Y}ang-{B}axter equation and invariants of
  links}, Invent. Math. 92 (1988) 527--553\relax \MR{939474}

\bibitem{Vi}
\textbf{O Viro}, \emph{Remarks on definition of Khovanov homology}, e-print (2002)
\arxiv{math.GT/0202199}

\bibitem{Yo}
\textbf{Y Yokota}, \emph{Skeins and quantum {${\rm SU}(N)$} invariants
  of {$3$}-manifolds}, Math. Ann. 307 (1997) 109--138\relax \MR{1427678}

\end{thebibliography}
\end{document}